\documentclass[11pt]{article}
\usepackage{geometry}                
\usepackage{graphicx}
\usepackage{amssymb}
\usepackage{epstopdf}
\usepackage{amsmath}
\usepackage{amstext}
\usepackage{amssymb}
\usepackage{amsthm}
\usepackage{txfonts}
\usepackage{eufrak}
\usepackage{dsfont}
\usepackage{yfonts}

\input xy
\xyoption{all}

\theoremstyle{plain}
\newtheorem{theorem}{Theorem}[section]
\newtheorem{definition}[theorem]{\bf Definition}
\newtheorem{proposition}[theorem]{Proposition}
\newtheorem{conjecture}[theorem]{Conjecture}
\newtheorem{remark}[theorem]{Remark}
\newtheorem{lemma}[theorem]{\bf Lemma}
\newtheorem{corollary}[theorem]{\bf Corollary}

\newcommand{\lm}{\Lambda}
\newcommand{\ra}{\rightarrow}

\newcommand{\mQ}{\mathbb{Q}}
\newcommand{\mZ}{\mathbb{Z}}

\newcommand{\mH}{\mathcal{H}} 
\newcommand{\mO}{\mathcal{O}}
\newcommand{\mC}{\mathcal{C}}
\newcommand{\G}{\Gamma}
\newcommand{\g}{\gamma}
\newcommand{\Gi}{\G^{(i)}}
\newcommand{\Ge}{\G^{(e)}}
\newcommand{\wh}{\widehat}
\newcommand{\ts}{\mathfrak{S}}
\newcommand{\R}{\mathfrak{R}}
\newcommand{\vp}{\varphi}
\newcommand{\N}{\mathcal{N}}
\newcommand{\mb}{\mathfrak{b}}

\newcommand{\ma}{\mathfrak{a}}
\newcommand{\mc}{\mathfrak{c}}

\newcommand{\ilim}[1]{%
	\displaystyle{%
	\lim_{\genfrac{}{}{0pt}{}{\longleftarrow}{\scriptstyle #1}} }\;}

\title{\sc Proof of the Main Conjecture of Noncommutative Iwasawa Theory for Totally Real Number Fields in Certain Cases}
\author{Mahesh Kakde}

\begin{document}
\maketitle{}
\begin{abstract} Fix an odd prime $p$. Let $G$ be a compact $p$-adic Lie group containing a closed, normal, pro-$p$ subgroup $H$ which is abelian and such that $G/H$ is isomorphic to the additive group of $p$-adic integers $\mZ_p$. First we assume that $H$ is finite and compute the Whitehead group of the Iwasawa algebra, $\lm(G)$, of $G$. We also prove some results about certain localisation of $\lm(G)$ needed in Iwasawa theory. Let $F$ be a totally real number field and let $F_{\infty}$ be an admissible $p$-adic Lie extension of $F$ with Galois group $G$. The computation of the Whitehead groups are used to show that the Main Conjecture for the extension $F_{\infty}/F$ can be deduced from certain congruences between abelian $p$-adic zeta functions of Delige and Ribet. We prove these congruences with certain assumptions on $G$. This gives a proof of the Main Conjecture in many interesting cases such as $\mZ_p\rtimes\mZ_p$-extensions.     
\end{abstract}
\tableofcontents
\section{Introduction}
Iwasawa theory studies the mysterious relationship between purely arithmetic objects and special values of complex $L$-functions. A precise form of this relationship is usually called the ``Main Conjecture". Historically, the first version of this Main Conjecture was stated for the cyclotomic $\mZ_p$-extension of a totally real number field. The aim of this paper is to prove the first noncommutative generalisation of this Main Conjecture.

I would like to express my gratitude to my advisor Professor John Coates for introducing me to this subject, for his continual encouragement, and many very helpful discussions and suggestions. I also thank Professor Kato for very generously sharing his insight of the subject during several valuable discussions, and also for providing me with a copy of \cite{Kato:2006}. I thank Takako Fukaya for providing me with the preprint of \cite{FukayaKato:2006}.  

\subsection{Iwasawa theory for totally real number fields}
In this section we review the Iwasawa theory for totally real number fields. Let $F$ be a totally real number field and $p$ be a prime. The field $F(\mu_{p^{\infty}}) = \cup_{n \geq 1} F(\mu_{p^n})$, where $ \mu_{p^n}$ denotes the group of $p^n$th roots of unity, contains a unique Galois extension of $F$ whose Galois group over $F$ is isomorphic to the additive group of $p$-adic integers $\mZ_p$. This extension is called the cyclotomic $\mZ_p$-extension of $F$ and we denote it by $F^{cyc}$. We write $\Gamma$ for the Galois group of $F^{cyc}$ over $F$, and fix a topological generator $\gamma$ of $\Gamma$. 

\begin{definition} An admissible $p$-adic Lie extension $F_{\infty}$ of $F$ is a Galois extension $F_{\infty}$ of $F$ such that ($i$) $F_{\infty}/F$ is unramified outside a finite set of primes of $F$, ($ii$) $F_{\infty}$ is totally real, ($iii$) $G:= Gal(F_{\infty}/F)$ is a $p$-adic Lie group, and ($iv$) $F_{\infty}$ contains $F^{cyc}$.
\label{admissible}
\end{definition}

\noindent Note that if $G$ has dimension $\geq 2$, Leopoldt's conjecture implies that $G$ must be non-abelian. Here is a typical example of such an admissible $p$-adic Lie extension with $G$ non-abelian. Let $F$ be the maximal real subfield of $\mQ(\mu_{37})$. Define $F_{\infty}$ to be the maximal abelian $37$-extension of $F^{cyc}$ which is unramified outside the unique prime above 37. Plainly, $F_{\infty}$ is Galois over $F$, and $G=Gal(F_{\infty}/F)$ has dimension 2 by classical computations on cyclotomic fields. 

\noindent From now on, $F_{\infty}/F$ will denote an admissible $p$-adic Lie extension, and we put
\[
G= Gal(F_{\infty}/F), \hspace{1cm} H=Gal(F_{\infty}/F^{cyc}), \hspace{1cm} \Gamma = G/H.
\]
Let $\Sigma$ denote any finite set of finite primes of $F$ which will always be assumed to contain the primes of $F$ which ramify in $F_{\infty}$.  \\

\noindent Throughout this paper $p$ will denote a fixed odd prime number. For any profinite group $P$, and $\mO$ the ring of integers of a finite extension of $\mQ_p$, we define
\[
\lm_{\mO}(P) = \ilim{U} \mO[P/U],
\]
where $U$ runs through the open normal subgroups of $P$, and $\mO[P/U]$ denotes the group ring of $P/U$ with coefficients in $\mO$. When $\mO = \mZ_p$, we write simply $\lm(P)$. Unless stated otherwise, we shall consider left modules over the Iwasawa algebras $\lm_{\mO}(P)$. Following Coates et.al. \cite{CFKSV:2005}, we define 
\[
S= \{s \in \lm(G) : \lm(G)/\lm(G)s \ \text{is a finitely generated} \ \lm(H)-\text{module} \}.
\]
It is proven in \cite{CFKSV:2005}, that $S$ is a multiplicatively closed subset of nonzero divisors in $\lm(G)$, which is left and right Ore set. Hence the localisation $\lm(G)_S$ of $\lm(G)$ exists, and the natural map from $\lm(G)$ to $\lm(G)_S$ is injective. A $\lm(G)$-module $M$ is called \emph{S-torsion} if every element of $M$ is annihilated by some element in $S$. It is proven in \emph{loc. cit.} that a $\lm(G)$-module $M$ is $S$-torsion if and only if it is finitely generated as a $\lm(H)$-module. \\

\noindent We now recall a few notions from algebraic $K$-theory. Most of this part is based on the introductory section of Fukaya-Kato \cite{FukayaKato:2006}. 
\begin{definition}For any ring $\lm$, $K_0(\lm)$ is an abelian group, whose group law we denote additively, defined by the following set of generators and relations. 
Generators : $[P]$, where $P$ is a finitely generated projective $\lm$-module. 
Relations: ($i$) $[P] = [Q]$ if $P$ is isomorphic to $Q$ as $\lm$ module, and ($ii$) $[P\oplus Q] = [P] + [Q]$.
\label{k0}
\end{definition}

\noindent It is easily seen that any element of $K_0(\lm)$ can be written as $[P]-[Q]$ for finitely generated projective $\lm$-modules $P$ and $Q$, and $[P]-[Q] = [P^{\prime}] - [Q^{\prime}]$ in $K_0(\lm)$ if and only if there is a finitely generated projective $\lm$-module $R$ such that $P\oplus Q^{\prime} \oplus R$ is isomorphic to $P^{\prime} \oplus Q \oplus R$. 

\begin{definition} For any ring $\lm$, $K_1(\lm)$ is an abelian group, whose group law we denote multiplicatively, defined by the following generators and relations. 
Generators : $[P, \alpha]$, where $P$ is a finitely generated projective $\lm$-module and $\alpha$ is an automorphism of $P$. 
Relations : ($i$) $[P, \alpha] = [Q, \beta]$ if there is an isomorphism $f$ from $P$ to $Q$ such that $f \circ \alpha = \beta \circ f$, ($ii$) $[P, \alpha \circ \beta] = [P, \alpha] [Q, \beta]$, and ($iii$) $[P\oplus Q, \alpha \oplus \beta] = [P, \alpha] [Q, \beta]$.
\label{k1}
\end{definition}

\noindent Here is an alternate description of $K_1(\lm)$. We have a canonical homomorphism $GL_n(\lm) \ra K_1(\lm)$ defined by mapping $\alpha$ in $GL_n(\lm)$ to $[\lm^n, \alpha]$, where $\lm^n$ is regarded as a set of row vectors and $\alpha$ acts on them from the right. Now using the inclusion maps $GL_n(\lm) \hookrightarrow GL_{n+1}(\lm)$ given by $g \mapsto \left( \begin{array}{cc} g & 0 \\ 0 & 1\end{array} \right)$, we let 
\[
GL_{\infty}(\lm) = \cup_{n \geq 1} GL_n(\lm).
\]
Then the canonical homomorphisms $GL_n(\lm) \ra K_1(\lm)$ induce an isomorphism (see for example \cite{Oliver:1988}, chapter 1)
\[
\frac{GL_{\infty}(\lm)}{[GL_{\infty}(\lm), GL_{\infty}(\lm)]} \xrightarrow{\sim} K_1(\lm),
\]
where $[GL_{\infty}(\lm), GL_{\infty}(\lm)]$ is the commutator subgroup of $GL_{\infty}(\lm)$. If $\lm$ is commutative, then the determinant maps, $GL_n(\lm) \ra \lm^{\times}$, induce the determinant map
\[
det: K_1(\lm) \ra \lm^{\times}, 
\]
via the above isomorphism. This gives a splitting of the canonical homomorphism $\lm^{\times} = GL_1(\lm) \ra K_1(\lm)$. If $\lm$ is semilocal then Vaserstein (\cite{Vaserstein:1969} and \cite{Vaserstein:2005}) proves that the canonical homomorphism $\lm^{\times} =GL_1(\lm) \ra K_1(\lm)$ is surjective. From these two facts we conclude that if $\lm$ is a semilocal commutative ring, then the determinant map induces a group isomorphism between $K_1(\lm)$ and $\lm^{\times}$. \\

\noindent Let $f: \lm \ra \lm^{\prime}$ be a ring homomorphism. We consider the category $\mC_f$ of all triplets $(P,\alpha, Q)$, where $P$ and $Q$ are finitely generated projective $\lm$-modules and $\alpha$ is an isomorphism between $\lm^{\prime} \otimes_{\lm}P$ and $\lm^{\prime} \otimes_{\lm} Q$ as $\lm^{\prime}$-modules. A morphism between $(P, \alpha, Q)$ and $(P^{\prime}, \alpha^{\prime}, Q^{\prime})$ is a pair of $\lm$-module morphisms $g: P \ra P^{\prime}$ and $h: Q \ra Q^{\prime}$ such that 
\[
\alpha^{\prime}\circ (Id_{\lm^{\prime}}\otimes g) = (Id_{\lm^{\prime}}\otimes h)\circ \alpha .
\]
It is an isomorphism if both $g$ and $h$ are isomorphisms. A sequence of maps
\[
0 \ra (P^{\prime}, \alpha^{\prime}, Q^{\prime}) \ra (P, \alpha, Q) \ra (P^{\prime\prime}, \alpha^{\prime\prime}, Q^{\prime\prime}) \ra 0,
\]
is a short exact sequence if the underlying sequences
\[
0 \ra P^{\prime} \ra P \ra P^{\prime\prime} \ra 0, \ \text{and} 
\]
\[
0 \ra Q^{\prime} \ra Q \ra Q^{\prime\prime} \ra 0,
\]
are short exact sequences.

\begin{definition} For any ring homomorphism $f: \lm \ra \lm^{\prime}$, $K_0(f)$ is an abelian group, whose group law we denote additively, defined by the following generators and relations. Generators : $[(P, \alpha, Q)]$, where $(P, \alpha, Q)$ is an object in $\mC_f$. Relations : ($i$) $[(P, \alpha, Q)] = [(P^{\prime}, \alpha^{\prime}, Q^{\prime})]$ if $(P, \alpha, Q)$ is isomorphic to $(P^{\prime}, \alpha^{\prime}, Q^{\prime})$, ($ii$) $[(P, \alpha, Q)] = [(P^{\prime}, \alpha^{\prime}, Q^{\prime})] + [(P^{\prime\prime}, \alpha^{\prime\prime}, Q^{\prime\prime})]$ for every short exact sequence as above, and ($iii$) $[(P_1, \beta\circ\alpha, P_3)] = [(P_1, \alpha, P_2)] + [(P_2, \beta, P_3)]$.
\label{relativek0}
\end{definition} 

\noindent For the canonical injection, say $i$, of $\lm(G)$ in $\lm(G)_S$, we write $K_0(\lm(G), \lm(G)_S)$ for $K_0(i)$ and call it the \emph{relative} $K_0$. We give two other descriptions of this group $K_0(\lm(G), \lm(G)_S)$. For details see Weibel \cite{Weibel:2000}. Let $\mC_S$ be the category of bounded complexes of finitely generated projective $\lm(G)$ modules whose cohomologies are $S$-torsion. Consider the abelian group $K_0(\mC_S)$ with the following set of generators and relations. Generators : $[C]$, where $C$ is an object of $\mC_S$. Relations : ($i$) $[C]=0$ if $C$ is acyclic, and ($ii$) $[C] = [C^{\prime}] + [C^{\prime\prime}]$, for every short exact sequence 
\[
0 \ra C^{\prime} \ra C \ra C^{\prime\prime} \ra 0, 
\]
in $\mC_S$. \\

\noindent Let $\mH_S$ be the category of finitely generated $\lm(G)$-modules which are $S$-torsion and which have a finite resolution by finitely generated projective $\lm(G)$-modules. Let $K_0(\mH_S)$ be the abelian group defined by the following set of generators and relations. Generators : $[M]$, where $M$ is an object of $\mH_S$. Relations : ($i$) $[M] =[M^{\prime}]$ if $M$ is isomorphic to $M^{\prime}$, and ($ii$) $[M] =[M^{\prime}] + [M^{\prime\prime}]$ for every short exact sequence 
\[
0 \ra M^{\prime} \ra M \ra M^{\prime\prime} \ra 0, 
\]
of modules in $\mH_S$. There are isomorphisms
\[
K_0(\lm(G), \lm(G)_S) \xrightarrow{\sim} K_0(\mC_S), \ \text{and}
\]
\[
K_0(\lm(G), \lm(G)_S) \xrightarrow{\sim} K_0(\mH_S),
\]
given as follows. First observe that every isomorphism $\alpha$ from $\lm(G)_S\otimes_{\lm(G)}P$ to $\lm(G)_S\otimes_{\lm(G)}Q$ is of the form $s^{-1}a$ with $a$ a $\lm(G)$-homomorphism from $P$ to $Q$ and $s$ an element of $S$. Then the above mentioned isomorphisms are respectively given by
\[
[(P, \alpha, Q)] \mapsto [P \xrightarrow{a} Q] + [Q \xrightarrow{s} Q], \ \text{and}
\]
\[
[(P, \alpha, Q)] \mapsto [Q/a(P)] + [Q/Qs].
\]
The relative $K_0$ fits into the localisation exact sequence of $K$-theory
\[
K_1(\lm(G)) \ra K_1(\lm(G)_S) \xrightarrow{\partial} K_0(\lm(G), \lm(G)_S) \xrightarrow{\eta} K_0(\lm(G)) \ra K_0(\lm(G)_S).
\]
The homomorphism $\partial$ maps $\alpha$ in $K_1(\lm(G)_S)$ to $[(\lm(G)^n, \tilde{\alpha}, \lm(G)^n)]$, where $\tilde{\alpha}$ is any lift of $\alpha$ in $GL_{\infty}(\lm(G)_S)$ and $\tilde{\alpha}$ lies in $GL_n(\lm(G)_S)$. The homomorphism $\eta$ maps $[(P, \alpha, Q)]$ to $[P] -[Q]$ in $K_0(\lm(G))$. The following lemma is essentially proven in \cite{CFKSV:2005} under the assumption that $G$ has no element of order $p$. The same technique gives this more general result as we now show.

\begin{lemma} The homomorphism $\partial$ is surjective.
\label{localisationsequence}
\end{lemma}
\noindent{\bf Proof:} We will show that the homomorphism $\eta$ is 0. Fix a pro-$p$ open normal subgroup $P$ of $G$, and put $\Delta = G/P$. We write $\mathcal{V} = \mathcal{V}(\Delta)$ for the set of irreducible representations of the finite group $\Delta$ over $\bar{\mQ}_p$ and we take $L$ to be some fixed finite extension of $\mQ_p$ such that all representations in $\mathcal{V}$ can be realised over $L$. Recall the construction in \cite{CFKSV:2005} of the canonical homomorphism $\lambda$ from $K_0(\lm(G))$ to $\prod_{\rho \in \mathcal{V}} K_0(L)$. It is the composition $\lambda = \lambda_4 \circ \lambda_3 \circ \lambda_2 \circ \lambda_1$ of four natural maps $\lambda_i$ ($i = 1,2,3,4$) as follows
\[
\lambda_1 : K_0(\lm(G)) \ra K_0(\mZ_p[\Delta]),
\]
\[
\lambda_2 : K_0(\mZ_p[\Delta]) \ra K_0(\mQ_p[\Delta]),
\]
\[
\lambda_3 : K_0(\mQ_p[\Delta]) \ra K_0(L[\Delta]),
\]
\[
\lambda_4 : K_0(L[\Delta]) \xrightarrow{\sim} \prod_{\rho \in \mathcal{V}}K_0(M_{n_{\rho}}(L)) \xrightarrow{\sim} \prod_{\rho \in \mathcal{V}} K_0(L).
\]
It can be shown that $\lambda_1$, $\lambda_2$, and $\lambda_3$ are injective. Hence $\lambda$ is injective as well. We now recall an alternate description of $\lambda$ given in \emph{loc. cit.}. The augmentation map from $\lm_{\mO}(G)$ to $\mO$ induces a map from $K_0(\lm_{\mO}(G))$ to $K_0(\mO)$ which we denote by $\tau$. Let $U$ be a finitely generated $\lm_{\mO}(G)$-module having a finite resolution by finitely generated projective $\lm_{\mO}(G)$-modules. Then one can define the class of $U$, denoted by $[U]$, in $K_0(\lm_{\mO}(G))$. Since $\mO$ is a domain, $K_0(\mO)$ can be identified with the Grothendieck group of the category of all finitely generated $\mO$-modules. Then $\tau$ is explicitly given by 
\[
\tau([U]) = \sum_{ i \geq 0} (-1)^i [H_i(G, U)].
\]
$\tau$ factors through the map
\[
\epsilon : K_0(\lm_{\mO}(G)) \ra K_0(\lm_{\mO}(\Gamma)),
\]
given by the natural surjection of $\lm_{\mO}(G)$ on $\lm_{\mO}(\Gamma)$. Moreover, $\epsilon$ is given explicitly by 
\[
\epsilon([U]) = \sum_{ i \geq 0} (-1)^i[H_i(H, U)].
\]
Now take $\mO$ to be the ring of integers of $L$ and $j$ be the isomorphism from $K_0(\mO)$ to $K_0(L)$ induced by the inclusion of $\mO$ in $L$. Let $tw_{\rho}(U)$ be defined by $U \otimes_{\mZ_p} \mO^{n_{\rho}}$, for any $\rho$ in $\mathcal{V}$. It can be made into a left $G$ module by the diagonal action. This action extends to make $tw_{\rho}(U)$ a $\lm(G)$-module. It is proven in \cite{CFKSV:2005} that $tw_{\rho}(U)$ is a finitely generated $S$-torsion $\lm(G)$-module if $U$ is. 

\noindent We now finish the proof. Take $[(P, \alpha, Q)]$ in $K_0(\lm(G), \lm(G)_S)$. As remarked earlier, $\alpha$ is of the form $s^{-1}a$, with $a$ a $\lm(G)$-morphism from $P$ to  $Q$ and $s$ is an element of $S$. We will show that $[Q/a(P)]$ is 0 in $K_0(\lm(G))$. Put $M= Q/a(P)$. Then $H_i(H, tw_{\rho}(M))$ is a $\lm_{\mO}(\Gamma)$-torsion module for all $i \geq 0$ and thus its class in $K_0(\lm_{\mO}(\Gamma))$ vanishes. Hence it follows that $\epsilon([tw_{\rho}(M)]) = 0$, whence $\tau([tw_{\rho}(M)]) = 0$ for all $\rho$ in $\mathcal{V}$. But this implies that $\lambda([M])=0$ and so $[M]=0$ in $K_0(\lm(G))$. This completes the proof. \qed \\

\noindent We now explain the basic arithmetic objects attached to $F_{\infty}/F$ which will be needed to formulate the Main Conjecture.

\begin{definition} We say that $F_{\infty}/F$ satisfies the hypothesis $\mu=0$ if there exists a pro-$p$ open subgroup $H^{\prime}$ of $H$ such that the Galois group over $F_{\infty}^{H^{\prime}}$ of the maximal unramified abelian $p$-extension of $F_{\infty}^{H^{\prime}}$ is a finitely generated $\mZ_p$-module. 
\label{hypothesisonmu}
\end{definition}

\noindent Of course, this is a special case of a general conjecture of Iwasawa asserting that, for every finite extension $K$ of $\mQ$, the Galois group over $K^{cyc}$ of the maximal unramified abelian $p$-extension of $K^{cyc}$ is a finitely generated $\mZ_p$-module. When $K$ is an abelian extension of $\mQ$, this conjecture is proven by Ferrero-Washington \cite{FerreroWashington:1979}. \\

\noindent Let $M_{\infty}$ be the maximal abelian $p$-extension of $F_{\infty}$, which is unramified outside the set of primes above $\Sigma$. As usual $G$ acts on $Gal(M_{\infty}/F_{\infty})$ by $g\cdot x = \tilde{g}x\tilde{g}^{-1}$, where $g$ is in $G$, and $\tilde{g}$ is any lifting of $g$ to $Gal(M_{\infty}/F)$. This action extends to a left action of $\lm(G)$. 

\begin{lemma} $Gal(M_{\infty}/F_{\infty})$ is finitely generated over $\lm(H)$ if and only if $F_{\infty}/F$ satisfies the hypothesis $\mu = 0$.
\label{xistorsion}
\end{lemma}
\noindent{\bf Proof:} Put $X = Gal(M_{\infty}/F_{\infty})$, and let $H^{\prime}$ be any pro-$p$ open subgroup of $H$. Thus $\lm(H^{\prime})$ is a local ring. It follows from Nakayama's lemma that $X$ is finitely generated over $\lm(H)$ if and only if $X_{H^{\prime}}$ is a finitely generated $\mZ_p$-module. Let $F_{\Sigma}$ denote the maximal pro-$p$ extension of $F$ which is unramified outside $\Sigma$, and let $K_{\infty} = F_{\infty}^{H^{\prime}}$. Then we have the inflation-restriction exact sequence
\[
0 \ra H^1(H^{\prime}, \mQ_p/\mZ_p) \ra H^1(Gal(F_{\Sigma}/K_{\infty}), \mQ_p/\mZ_p) \ra H^1(Gal(F_{\Sigma}/F_{\infty}), \mQ_p/\mZ_p)^{H^{\prime}} \ra H^2(H^{\prime}, \mQ_p/\mZ_p).
\]
As $H^{\prime}$ is a $p$-adic Lie group, $H^i(H^{\prime}, \mQ_p/\mZ_p)$ are cofinitely generated $\mZ_p$-modules for all $i \geq 0$. Moreover, since $Gal(F_{\Sigma}/F)$ acts trivially on $\mQ_p/\mZ_p$, we have 
\[
H^1(Gal(F_{\Sigma}/L), \mQ_p/\mZ_p) = Hom(Gal(M_L/L), \mQ_p/\mZ_p),
\]
for every intermediate field $L$ with $F_{\Sigma} \supset L \supset F$; here $M_L$ denotes the maximal abelian $p$-extension of $L$ which is unramified outside $\Sigma$. We conclude from the above exact sequence that $X_{H^{\prime}}$ is a finitely generated $\mZ_p$-module if and only if $Gal(M_{K_{\infty}}/K_{\infty})$ is a finitely generated $\mZ_p$-module. Let $J_{\infty}$ denote the maximal unramified abelian $p$-extension of $K_{\infty}$. Now $K_{\infty}$ is clearly the cyclotomic $\mZ_p$-extension of some totally real finite extension $K$ of $F$, and Iwasawa \cite{Iwasawa:1973} has proven the weak Leopoldt conjecture holds for $K_{\infty}/K$. It follows easily from the validity of the weak Leopoldt conjecture for the totally real extension $K_{\infty}/K$ that $Gal(M_{\infty}/J_{\infty})$ is always a finitely generated $\mZ_p$-module. The conclusion of the lemma is now plain since $Gal(J_{\infty}/K_{\infty})$ being a finitely generated $\mZ_p$-module is precisely the hypothesis $\mu = 0$. \qed \\

\noindent In the view of the lemma, we shall henceforth make the following \\
{\bf Assumption:} Our admissible $p$-adic Lie extension $F_{\infty}/F$ satisfies the hypothesis $\mu = 0$. \\

\noindent $X =Gal(M_{\infty}/F_{\infty})$ is a fundamental arithmetic object which is studied through the Main Conjecture. The above lemma shows that $X$ is a $S$-torsion module, however, if $G$ has elements of order $p$, then $X$ may not have a finite resolution by finitely generated projective $\lm(G)$-modules. In our approach to the proof of the Main Conjecture it is necessary to consider $G$ having elements of order $p$. So we consider a complex of $\lm(G)$-modules which is closely related to $X$, and is quasi-isomorphic to a bounded complex of finitely generated projective $\lm(G)$-modules. We put $C(F_{\infty}/F)$ to be the complex
\[
C(F_{\infty}/F) = Hom(R\G_{\acute{e}t}(Spec(\mO_{F_{\infty}}[\frac{1}{\Sigma}]), \mQ_p/\mZ_p), \mQ_p/\mZ_p).
\]
Here $\mQ_p/\mZ_p$ is the locally constant sheaf corresponding to the abelian group $\mQ_p/\mZ_p$ on the \'{e}tale site of $Spec(\mO_{F_{\infty}}[\frac{1}{\Sigma}])$. Since $\mQ_p/\mZ_p$ is a direct limit of finite abelian groups of $p$-power order, we have an isomorphism
\[
R\G_{\acute{e}t}(Spec(\mO_{F_{\infty}}[\frac{1}{\Sigma}]), \mQ_p.\mZ_p) \xrightarrow{\sim} R\G(Gal(F_{\Sigma}/F_{\infty}), \mQ_p/\mZ_p).
\]
Recall that $F_{\Sigma}$ is the maximal $p$-extension of $F$ unramified outside $\Sigma$. Then $H^i(C(F_{\infty}/F))$ is 0 unless $ i $ is 0 or -1. $H^0(C(F_{\infty}/F)) = \mZ_p$ and $H^{-1}(C(F_{\infty}/F)) =Gal(M_{\infty}/F_{\infty})$. The following results are proven in Fukaya-Kato \cite{FukayaKato:2006}. \\
\noindent ($i$) The complex $C(F_{\infty}/F)$ is quasi-isomorphic to a bounded complex of finitely generated projective $\lm(G)$-modules. By passing to the derived category we identify $C(F_{\infty}/F)$ with a quasi-isomorphic bounded complex of finitely generated $\lm(G)$-modules. \\
\noindent ($ii$) If $F \subset K \subset F_{\infty}$ is any extension of $F$, then 
\[
\lm(Gal(K/F)) \otimes_{\lm(G)}C(F_{\infty}/F) \xrightarrow{\sim} C(K/F), 
\]
where $\lm(G)$ acts on the right on $\lm(Gal(K/F))$ through the natural surjection of $\lm(G)$ on $\lm(Gal(K/F))$. \\

\noindent By lemma \ref{xistorsion}, we know that $C(F_{\infty}/F)$ is $S$-torsion i.e. $\lm(G)_S \otimes_{\lm(G)} C(F_{\infty}/F)$ is acyclic. Hence we can talk about the class of $C(F_{\infty}/F)$, $[C(F_{\infty}/F)]$, in $K_0(\lm(G), \lm(G)_S)$. \\

\noindent We now explain the analytic objects in the Iwasawa theory of $F_{\infty}/F$ which we need to formulate the Main Conjecture. 

\noindent Let $\rho$ be an Artin representation (i.e. kernel of $\rho$ is open) of $Gal(\bar{F}/F)$ on a finite dimensional vector space over $\bar{\mQ}_p$, factoring through $G$. Let $\alpha$ be an embedding of $\bar{\mQ}_p$ in $\mathbb{C}$. The resulting complex Artin representation $\alpha \circ \rho$ gives the complex Artin $L$-function $L(\alpha \circ \rho, s)$. A famous result of Siegel says that $L(\alpha \circ \rho, -n)$ is an algebraic number for any odd positive integer $n$. This number depends on the choice of $\alpha$. If $\alpha$ is replaced by another embedding then we get a conjugate of the algebraic number $L(\alpha \circ \rho, -n)$ over $\mQ$. However, the beauty of the result of Siegel is that it makes it possible to choose a canonical conjugate. Hence we obtain an algebraic number $L(\rho, -n)$, \emph{``the value of complex $L$-function associated to $\rho$ at $-n$"}. For details see section 1.2 in Coates-Lichtenbaum \cite{CoatesLichtenbaum:1973}. Similarly, we can define the value of complex $L$-function associated to $\rho$ at $-n$ with Euler factors at primes in $\Sigma$ removed. We denote it by $L_{\Sigma}(\rho, -n)$. This value obviously depends only on the character associated to $\rho$. We will use the same letter to denote Artin representation and its character.  \\

\noindent Let $L$ be a finite extension of $\mQ_p$ with ring of integers $\mO$. Let $\rho$ be a continuous homomorphism from $G$ into $GL_n(\mO)$. It induces a ring homomorphism from $\lm(G)$ into $M_n(\lm_{\mO}(\G))$. This homomorphism is given on elements of $G$ by mapping $\sigma$ in $G$ to $\rho(\sigma)\bar{\sigma}$, where $\bar{\sigma}$ is the image of $\sigma$ in $\Gamma$. We write $Q_{\mO}(\G)$ for the field of fraction of $\lm_{\mO}(\G)$. It is proven in \cite{CFKSV:2005} that this homomorphism extends to a homomorphism 
\[
\Phi_{\rho} : \lm(G)_S \ra M_n(Q_{\mO}(\G)).
\]
$\Phi_{\rho}$ induces a homomorphism 
\[
\Phi_{\rho}^{\prime} : K_1(\lm(G)_S) \ra K_1(M_n(Q_{\mO}(\G))) = Q_{\mO}(\G)^{\times}.
\]
Now let $\varphi$ be the augmentation map from $\lm_{\mO}(\G)$ to $\mO$. We denote its kernel by $\mathfrak{p}$. If we write $\lm_{\mO}(\G)_{\mathfrak{p}}$ for the localisation of $\lm_{\mO}(\G)$ at the prime ideal $\mathfrak{p}$, then $\varphi$ extends to a homomorphism 
\[
\varphi : \lm_{\mO}(\G)_{\mathfrak{p}} \ra L.
\]
We extend this map to a map $\varphi^{\prime}$ from $Q_{\mO}(\G)$ to $L \cup \{\infty\}$ by mapping any $x$ in $Q_{\mO}(\G) - \lm_{\mO}(\G)_{\mathfrak{p}}$ to $\infty$. The composition of $\Phi_{\rho}^{\prime}$ with $\varphi^{\prime}$ gives a map
\begin{center}
$K_1(\lm(G)_S) \ra L \cup \{\infty\}$ \\
$ x \mapsto x(\rho) $.
\end{center}
This map has the following properties: \\
($i$) Let $G^{\prime}$ be an open subgroup of $G$. Let $\chi$ be an one dimensional representation of $G^{\prime}$ and $\rho = Ind_{G^{\prime}}^{G}(\chi)$. If $N$ is the norm map from $K_1(\lm(G)_S)$ to $K_1(\lm(G^{\prime})_S)$, then for any $x$ in $K_1(\lm(G)_S)$, $x(\rho) = N(x)(\chi)$. \\
($ii$) Let $\rho_i$ be continuous homomorphisms from $G$ into $GL_{n_i}(L)$, for $i=1,2$. We then get a continuous homomorphism $\rho_1\oplus \rho_2$ from $G$ into $GL_{n_1+n_2}(L)$. Then for any $x$ in $K_1(\lm(G)_S)$, $x(\rho_1\oplus \rho_2) = x(\rho_1)x(\rho_2)$. \\
($iii$) Let $U$ be a subgroup of $H$, normal in $G$. Let $\pi$ be the homomorphism from $K_1(\lm(G)_S)$ to $K_1(\lm(G/U)_S)$ induced by the natural surjection of $\lm(G)_S$ onto $\lm(G/U)_S$. Let $\rho$ be a continuous homomorphism of $G/U$ into $GL_n(L)$. We write $inf(\rho)$ for the composition of the natural surjection from $G$ onto $G/U$ with $\rho$. Then for any $x$ in $K_1(\lm(G)_S)$, $x(inf(\rho)) = \pi(x)(\rho)$. \\

\noindent{\bf Notation:} Let $\kappa_F$ be the cyclotomic character from $Gal(F(\mu_{p^{\infty}})/F)$ into $\mZ_p^{\times}$ given by 
\[
\sigma(\zeta) = \zeta^{\kappa_{F}(\sigma)},
\]
for any $\sigma$ in $Gal(F(\mu_{p^{\infty}})/F)$ and any $p$ power root of unity $\zeta$. 

\begin{conjecture} There exists $\zeta(F_{\infty}/F)$ in $K_1(\lm(G)_S)$ such that for any Artin character $\rho$ of $G$ and any positive integer $r$ divisible by $p-1$,
\[
\zeta(F_{\infty}/F)(\rho\kappa_F^r) = L_{\Sigma}(\rho, 1-r).
\]
\label{zeta}
\end{conjecture}

\begin{remark} $\zeta(F_{\infty}/F)$ is called a $p$-adic zeta function for the extension $F_{\infty}/F$. It depends on $\Sigma$ but we suppress this fact in the notation. It is also conjectured that $p$-adic zeta function is unique if it exists. If $G$ is abelian then existence and uniqueness of the $p$-adic zeta function is well known, as we will soon see.
\end{remark}

\noindent We are now ready to state the Main Conjecture. Recall the localisation sequence of $K$-theory
\[
K_1(\lm(G)) \ra K_1(\lm(G)_S) \xrightarrow{\partial} K_0(\lm(G), \lm(G)_S) \ra 0.
\]

\begin{conjecture} ({\bf Main Conjecture}) There is a unique element $\zeta(F_{\infty}/F)$ in $K_1(\lm(G)_S)$ such that $\partial(\zeta(F_{\infty}/F)) = -[C(F_{\infty}/F)]$ and for any Artin character $\rho$ of $G$ and any positive integer $r$ divisible by $p-1$,
\[
\zeta(F_{\infty}/F)(\rho\kappa_F^r) = L_{\Sigma}(\rho, 1-r).
\]
\label{mainconjecture}
\end{conjecture}

\begin{remark} The Main Conjecture in this form was first formulated by Kato \cite{Kato:1993}, and Fontaine and Perrin-Riou \cite{FontaineRiou:1991} in the case when $G$ is abelian. This was generalised to include noncommutative groups $G$ by Burns and Flach \cite{BurnsFlach:2001}, Huber and Kings \cite{HuberKings:2002}, Coates et.al. \cite{CFKSV:2005}, and Fukaya and Kato \cite{FukayaKato:2006}. Ritter and Weiss \cite{RitterWeiss:2} considered the case of one dimensional $p$-adic Lie groups. 
\end{remark}
\begin{remark} Of course, the Main Conjecture contains conjecture \ref{zeta}. 
\end{remark}
\begin{remark} One can show that the validity of the Main Conjecture is independent of $\Sigma$ as long as it contains all primes of $F$ which ramify in $F_{\infty}$. 
\end{remark}

\subsection{The commutative case} We now assume that $G$ is abelian. In this case conjectures \ref{zeta} and \ref{mainconjecture} are known to be true thanks to the deep works of Kubota-Leopoldt, Iwasawa, Deligne-Ribet, Mazur-Wiles, and Wiles, among others. Though all the results in this section are well-known, we collect them as we need them later. First we introduce some notations. As $G$ is abelian, we have $G \cong H\times \G$. We also assume that $H$ is a finite group. Of course, if the Leopoldt conjecture is true for $F$ and $p$, then $H$ has to be finite. The isomorphism $G \cong H \times\G$ gives an isomorphism of $\lm_{\mO}(G)$ with $\cong \mO[H][[T]]$, the power series ring in the indeterminate $T$ with coefficients in $\mO[H]$, obtained by mapping $\g$ to $T+1$ (recall that $\g$ is the fixed topological generator of $\G$). For a pro-finite group $P$, we put $\hat{P}$ for the set of all one dimensional $p$-adic character of $P$ of finite order. For any $\psi \in \hat{H}$, we put $\mO_{\psi}$ for the ring of integers in the field obtained by adjoining all the values of $\psi$ to $\mQ_p$. We consider $\psi$ as a character of $G$ by using the surjection of $G$ onto $H$. We write $\tilde{\psi}$ for the homomorphism 
\[
\lm(G) \ra \lm_{\mO_{\psi}}(\G) \cong \mO_{\psi}[[T]],
\]
induced by $\sigma \mapsto \psi(\sigma)\bar{\sigma}$, for every $\sigma$ in $G$. Here $\bar{\sigma}$ denotes the image of $\sigma$ in $\G$. If $r$ is a positive integer divisible by $p-1$, then $\kappa_{F}^r$ is a continuous homomorphism from $\G$ to $\mZ_p^{\times}$. For any $f \in \lm_{\mO_{\psi}}(\G)$, we can evaluate $f$ at $\kappa_F^{r}$ to get $f(\kappa_{F}^r)$. We remark that $f(\kappa_F^{r})$ is traditionally written as $\int_{\G} \kappa_F^rdf$ or $\kappa_F^r(f)$ in commutative Iwasawa theory (see for example Coates-Sujatha \cite{CoatesSujatha:2007}). We are following the notation of \cite{CFKSV:2005}. Under the isomorphism $\lm_{\mO_{\psi}}(\G) \xrightarrow{\vp} \mO_{\psi}[[T]]$, the compatibility of evaluation at $\kappa_F^r$ is given by 
\[
f(\kappa_F^r) = \vp(f)(\kappa_F^r(\g)-1).
\]
The following statements are known to be true: \\
($i$) Let $\psi$ be an element of $\hat{G}$. Let $F_{\psi}$ be the fixed field of the kernel of $\psi$. Following Greenberg, $\psi$ is called of type $W$ if $F_{\psi} \subset F^{cyc}$. Then there exists a power series $G_{\psi, \Sigma}(T)$ in $O_{\psi}[[T]]$ such that for every positive integer $r$  divisible by $p-1$, we have
\[
\frac{G_{\psi, \Sigma}(\kappa_F^r(\g)-1)}{H_{\psi}(\kappa_F^r(\g)-1)} = L_{\Sigma}(\psi, 1-r),
\]
where $H_{\psi}(T)$ is the polynomial $\psi(\g)(T+1)-1$ if $\psi$ is of type $W$ and is $1$ otherwise. 

\noindent For $F=\mQ$, this was first proven by Iwasawa, following the work of Kubota-Leopoldt. For a totally real quadratic extension $F$ of $\mQ$, this was proven by Coates-Sinnott \cite{CoatesSinnott:1974} using the explicit formulae of Siegel for values of partial zeta functions at non-positive integers. Coates \cite{Coates:1977} gave certain hypothesis, in terms of congruences between the values of partial zeta functions, which gives $G_{\psi, \Sigma}(T)$ in general. These hypothesis were proven by Deligne-Ribet \cite{DeligneRibet:1980}, and also by Cassou-Nogu\'{e}s \cite{Cassou:1979}. Deligne and Ribet used $p$-adic Hilbert modular forms to construct $G_{\psi, \Sigma}(T)$ generalising the work of Serre \cite{Serre:1973} who constructed $G_{\psi, \Sigma}(T)$ for the case when $\psi$ is a power of \emph{Teichm\"{u}ller} character. Cassou-Nogu\'{e}s used explicit formulae of Shintani \cite{Shintani:1976} to prove the hypothesis of Coates.  

\noindent An element $f$ in the field of fractions of $\lm(G)$ is called a \emph{pseudomeasure} if $(g-1)f$ lies in $\lm(G)$ for every $g$ in $G$. Serre's account \cite{Serre:1978} of the work of Deligne and Ribet shows that there is a unique pseudomeasure $\zeta(F_{\infty}/F)$ such that for every $\psi \in \hat{G}$, we have 
\[
\tilde{\psi}(\zeta(F_{\infty}/F)) = \frac{G_{\psi, \Sigma}(T)}{H_{\psi}(T)}.
\]

\noindent ($ii$) Let $\psi \in \hat{H}$. Put $V= \bar{\mQ}_p \times_{\mZ_p} X$. Then $V$ is a finite dimensional $\bar{\mQ}_p$ vector space (this follows from lemma \ref{xistorsion} but is known to be true, due to Iwasawa, even without the assumption $\mu = 0$ for $F_{\infty}/F$). Let $h_{\psi}(T)$ be the characteristic polynomial of $(\g-1)$ acting on the space $V^{\psi} = \{ v \in V : hv = \psi(h)v \text{   for all   } h \in H\}$. If $f(T)$ is any power series in $\mO_{\psi}[[T]]$, then by Weierstrass Preparation Theorem we can write $f(T)$ as $\pi^{\mu(f)}f^*(T)u(T)$, where $\pi$ is a uniformiser of $\mO_{\psi}$, $f^*(T)$ is a distinguished polynomial and $u(T)$ is a unit in $\mO_{\psi}[[T]]$. Then we have
\[
h_{\psi}(T) = G_{\psi, \Sigma}^*(T) \hspace{2cm} \text{(Iwasawa Main Conjecture)}.
\]
This is the Iwasawa Main Conjecture for totally real number fields. This was formulated for $F =\mQ$ by Iwasawa, and extended to arbitrary totally real number fields by Coates \cite{Coates:1977} and by Greenberg \cite{Greenberg:1983}. For $F=\mQ$ this conjecture is proven by Mazur-Wiles \cite{MazurWiles:1984} after deep results in this direction by Iwasawa. The general result is proven by Wiles \cite{Wiles:1990}. For $F=\mQ$, Rubin \cite{Lang:1990} gave another proof using Kolyvagin's technique of Euler systems. \\

\noindent We state an important result about the $\mu$ invariant. Let $\psi \in \hat{H}$ be of order prime to $p$. We write $X^{\psi}$ for 
\[
(X \otimes_{\mZ_p} \mO_{\psi})^{(\psi)} := \{ x \in X \otimes_{\mZ_p}\mO_{\psi} : \sigma x = \psi(\sigma) x \text{ for all } \sigma \in H \}.
\]
Then $X^{\psi}$ is a finitely generated torsion $\mO_{\psi}[[T]]$-module. Structure theory of such modules gives a characteristic power series which in this case has a form $\pi^{\mu(X^{\psi})}h_{\psi}(T)$. And Wiles (\cite{Wiles:1990}, th. 1.4) shows that 
\begin{equation}
\mu(X^{\psi}) = \mu(G_{\psi, \Sigma}(T)).
\label{muequation}
\end{equation}
Hence our assumption $\mu=0$ shows that both these quantities are 0. \\

\noindent These deep results immediately imply the Main Conjecture in the abelian case, as we sketch below. 

\begin{lemma} Let $F_{\infty}/F$ be an abelian admissible $p$-adic Lie extension satisfying the hypothesis $\mu = 0$. Then the $p$-adic zeta function of Deligne and Ribet, $\zeta(F_{\infty}/F)$, is a unit in $\lm(G)_S$. 
\label{abelianzetaisunit}
\end{lemma}
\noindent{\bf Proof:} We write $H$ as $H^{\prime} \times H_p$, where $H_p$ is the $p$-part of $H$ and $H^{\prime}$ is a finite group whose order is prime to $p$. Note that $\lm(G)_S$ in this case is just the localisation at the prime ideal $(p)$, generated by $p$. We have the following decomposition of $\lm(G)_{(p)}$
\[
\lm(G)_{(p)} \xrightarrow{\sim} \mZ_p[H^{\prime}\times H_p][[T]]_{(p)} \xrightarrow{\sim} \oplus_{\psi \in \hat{H^{\prime}}} \mO_{\psi}[H_p][[T]]_{\mathfrak{p}_{\psi}},
\]
where $\mathfrak{p}_{\psi}$ is the ideal in $\mO_{\psi}[H_p][[T]]$ generated by a uniformiser in $\mO_{\psi}$.  We must show that the image of $\zeta(F_{\infty}/F)$ in each summand is a unit. But each summand $\mO_{\psi}[H_p][[T]]_{\mathfrak{p}_{\psi}}$ is a local ring with maximal ideal 
\[
\mathfrak{m} = \{ x \in \mO_{\psi}[H_p][[T]]_{\mathfrak{p}_{\psi}} : \tilde{\mathds{1}}_{H_p}(x) \in \mathfrak{p}_{\psi} \},
\]
Hence we must show that 
\[
\zeta_{\psi} = (\tilde{\mathds{1}_{H_p}\times \psi})(\zeta(F_{\infty}/F)) \notin \mathfrak{p}_{\psi}.
\]
But $\zeta_{\psi}$ is none other than the power series $G_{\psi, \Sigma}(T)$ which does not lie in $\mathfrak{p}_{\psi}$ by the above remark about Wiles' result on the $\mu$ invariant (see equation \ref{muequation}). \qed \\

\noindent Consider the complex $C_{\psi} = \mO_{\psi}[[T]] \otimes_{\lm(G)} C(F_{\infty}/F)$. $C_{\psi}$ is a bounded complex of finitely generated projective $\mO_{\psi}[[T]]$-modules and is $(\mO_{\psi}[[T]] - \mathfrak{p}_{\psi})$-torsion. Hence we can talk about the class of $C_{\psi}$, $[C_{\psi}]$, in $K_0(\mO_{\psi}[[T]]. \mO_{\psi}[[T]]_{\mathfrak{p}_{\psi}})$. In fact $[C_{\psi}]$ is the image of $[C(F_{\infty}/F)]$ under the natural map
\[
K_0(\lm(G), \lm(G)_{(p)}) \ra K_0(\mO_{\psi}[[T]], \mO_{\psi}[[T]]_{\mathfrak{p}_{\psi}}).
\]
It is also easy to see that $H^i(C_{\psi})$ is 0 unless $i=0$ or $i=-1$ and $H^0(C_{\psi}) = \mO_{\psi}[[T]] \otimes_{\lm(G)} \mZ_p$ (which is 0 unless $\psi$ is trivial), and $H^{-1}(C_{\psi}) = \mO_{\psi}[[T]] \otimes_{\lm(G)}X$.  \\

\noindent We also make the following observation. If $Y$ is a finitely generated $\mO_{\psi}[[T]]$-module which is $(\mO_{\psi}[[T]] - \mathfrak{p}_{\psi})$-torsion, and if $f$ is a characteristic power series of $Y$, then $f$ lies in $\mO_{\psi}[[T]]_{\mathfrak{p}_{\psi}}^{\times}$. The class of $Y$ in $K_0(\mO_{\psi}[[T]], \mO_{\psi}[[T]]_{\mathfrak{p}_{\psi}})$ is given by $[(Y, 0 , \mO_{\psi}[[T]]/\mO_{\psi}[[T]]f)]$ and $f$ maps to the class of $Y$ under the connecting homomorphism 
\[
\partial : K_1(\mO_{\psi}[[T]]_{\mathfrak{p}_{\psi}}) \cong \mO_{\psi}[[T]]_{\mathfrak{p}_{\psi}}^{\times} \ra K_0(\mO_{\psi}[[T]], \mO_{\psi}[[T]]_{\mathfrak{p}_{\psi}}).
\]

\noindent The image of $\zeta(F_{\infty}/F)$ under the natural map from $\lm(G)_{(p)}^{\times}$ to $\mO_{\psi}[[T]]_{\mathfrak{p}_{\psi}}^{\times}$ is $\frac{G_{\psi, \Sigma}(T)}{H_{\psi}(T)}$. Notice that $H_{\psi}(T)$ is 1 unless $\psi$ is trivial (in which case it is $T$). Since $\mO_{\psi}[[T]]\otimes_{\lm(G)} X$ is a finitely generated $\mZ_p$-module (thanks to our assumption that $\mu = 0$ for $F_{\infty}/F$), $h_{\psi}(T)$ is a characteristic power series of $\mO_{\psi}[[T]]\otimes_{\lm(G)} X$. Using Wiles' theorem i.e. Iwasawa Main Conjecture, we conclude that 
\[
\partial(\tilde{\psi}(\zeta(F_{\infty}/F))) = -[C_{\psi}].
\]

\begin{theorem} (Wiles) Let $F_{\infty}/F$ be an abelian admissible $p$-adic Lie extension satisfying the hypothesis $\mu = 0$. Then under the connecting homomorphism 
\[
\partial: K_1(\lm(G)_{(p)})= \lm(G)_{(p)}^{\times} \ra K_0(\lm(G), \lm(G)_{(p)}),
\]
$\zeta(F_{\infty}/F)$ maps to $-[C(F_{\infty}/F)]$ i.e. the Main Conjecture holds in the case when $G$ is abelian. 
\label{abelianmainconjecture}
\end{theorem}
\noindent{\bf Proof:} Consider the following commutative diagram with exact rows and columns

\xymatrix{ & 0 \ar[d] & 0 \ar[d] &  &  \\
0 \ar[r] & \lm(G)^{\times} \ar[d]_{\theta} \ar[r] & \lm(G)_{(p)}^{\times} \ar[d]_{\theta_S} \ar[r]^{\partial} & K_0(\lm(G), \lm(G)_{(p)}) \ar[r] \ar[d]_{\theta_0} & 0 \\
0 \ar[r] & \prod \mO_{\psi}[[T]]^{\times} \ar[r] & \prod \mO_{\psi}[[T]]_{\mathfrak{p}_{\psi}}^{\times} \ar[r]^{\partial} & \prod K_0(\mO_{\psi}[[T]]. \mO_{\psi}[[T]]_{\mathfrak{p}_{\psi}}) \ar[r] & 0 } 

\noindent where the products range over all $\psi$ in $\hat{H}$. The morphisms $\theta$ and $\theta_S$ take $x$ to $(\tilde{\psi}(x))_{\psi}$. We now show that $\theta_0$ is injective. The image of $\theta$ (resp. $\theta_S$) consists of all tuples $(a_{\psi})$ in $\prod_{\psi} \mO_{\psi}[[T]]^{\times}$ (resp. $\prod_{\psi} \mO_{\psi}[[T]]_{\mathfrak{p}_{\psi}}^{\times}$) such that 
\[
\frac{1}{|H|} \sum_{h \in H}h \big( \sum_{\psi \in \hat{H}} a_{\psi} \psi(h^{-1}) \big)
\]
lies in $\mZ_p[H][[T]]^{\times}$ (resp. $\mZ_p[H][[T]]_{(p)}^{\times}$). Now is it easily seen that if an element $x$ in $\mZ_p[H][[T]]_{(p)}^{\times}$ has image in $\prod_{\psi \in \hat{H}} \mO_{\psi}[[T]]^{\times}$, then $x$ must actually lie in $\mZ_p[H][[T]]^{\times}$. Hence we conclude that 
\[
Image(\theta_S) \cap \prod_{\psi \in \hat{H}} \mO_{\psi}[[T]]^{\times} = Image(\theta).
\]
Snake lemma now gives that $\theta_0$ is injective. The result now follows because 
\[
\partial(\tilde{\psi}(\zeta(F_{\infty}/F))) = -[C_{\psi}] \ \text{  for all  } \psi \in \hat{H}, \text{  and}
\]
\[
\theta_0([C(F_{\infty}/F)]) = ([C_{\psi}])_{\psi}.
\]
\qed \\

\subsection{A proposed strategy for proving the Main Conjecture}
\subsubsection{$K_1$ groups}
It is evident that computation of $K_1$ groups of $\lm(G)$ and $\lm(G)_S$ is of great importance in Iwasawa theory. If $G$ is abelian, $K_1(\lm(G)) \cong \lm(G)^{\times}$ and $K_1(\lm(G)_S) \cong \lm(G)_S^{\times}$ and hence are well understood. In \cite{Kato:2005}, Kato suggested studying these $K_1$ groups for a noncommutative group $G$ through its ``abelian subquotients". More precisely, let $(U,V)$ be a pair such that $U$ is an open subgroup of $G$, $V$ is an open subgroup of $H$ and a normal subgroup of $U$ such that $U/V$ is abelian. Then we have a homomorphism
\[
\theta_{U,V} : K_1(\lm(G)) \ra K_1(\lm(U)) \ra K_1(\lm(U/V)) \cong \lm(U/V)^{\times},
\]
where the first map is the norm homomorphism and the second map is the one induced by the natural surjection of $\lm(U)$ on $\lm(U/V)$. Let $I$ be a set of such pairs $(U,V)$. Consider the homomorphism 
\[
\theta_{I} := (\theta_{U,V})_{(U,V) \in I}: K_1(\lm(G)) \ra \prod_{(U,V) \in I} \lm(U/V)^{\times}.
\]
By taking $I$ sufficiently large we may hope that $\theta_I$ is injective and at the same time describe its image. In \cite{Kato:2005}, Kato implemented this strategy to get a description for open subgroups of $\mZ_p\rtimes \mZ_p^{\times}$. The proof is rather technical and long. Later Kato himself outlined a much more elegant approach to this question using the integral logarithm of Oliver and Taylor (See Oliver, \cite{Oliver:1988}), and a student of Kato gave a generalisation of the result of \cite{Kato:2005} via this method. Kato (unpublished, \cite{Kato:2006}) also applied this technique to the case when  $G$ is a quotient of the $p$-adic Heisenberg group. We also note that homomorphisms analogous to $\theta_{U,V}$ and $\theta_{I}$ (which we denote by $\theta_{U,V,S}$ and $\theta_{I,S}$ respectively) can be defined on $K_1(\lm(G)_S)$ and some results about these maps are also established by Kato in \cite{Kato:2005} and \cite{Kato:2006}. In this paper we compute $K_1(\lm(G))$ and give some results about $K_1(\lm(G)_S)$ in the case when $G=H\rtimes \Gamma$, with $H$ a finite abelian $p$ group. 

\subsubsection{A key observation} Kato, following a beautiful observation of Burns, suggested the following strategy for proving the noncommutative Main Conjecture. This result uses the deep result of Wiles. 

\noindent We assume that there is a set $I$ of pairs $(U,V)$ as above such that:\\
C1) For any Artin character $\rho$ of $G$, there is a finite family $\{(U_i,V_i)\}$ in $I$ and one dimensional Artin character $\chi_i$ of $U_i/V_i$ such that $\rho$ is a $\mZ$-linear combination of $Ind_{U_i}^{G}\chi_i$. \\
C2) There is a subgroup $\Phi$ of $\prod_{(U,V)\in I}\lm(U/V)^{\times}$ such that $\theta$ induces an isomorphism 
\[
\theta: K_1(\lm(G)) \ra \Phi.
\]
C3) There is a subgroup $\Phi_S$ of $\prod_{(U,V)\in I}\lm(U/V)_S^{\times}$ such that $\Phi_S \cap \prod_{(U,V)\in I}\lm(U/V)^{\times} =\Phi$ and $\theta_S(K_1(\lm(G)_S)) \subset \Phi_S$. \\

\begin{proposition} (Burns, Kato) Let $F_{\infty}/F$ be an admissible $p$-adic Lie extension satisfying the hypothesis $\mu =0$. Suppose we have an $I$ which satisfies C1), C2) and C3) above. Then  the Main Conjecture is true for $F_{\infty}/F$ if and only if $(\zeta(F_{\infty}^V/F_{\infty}^U)) \in \Phi_S$.
\end{proposition}
\noindent{\bf Proof:} If the Main Conjecture is true then, in particular, $\zeta(F_{\infty}/F)$ exists and $\theta_S(\zeta(F_{\infty}/F)) = (\theta_{U,V,S}(\zeta(F_{\infty}/F))) \in \Phi_S$. $\theta_{U,V,S}(\zeta(F_{\infty}/F)) = \zeta(F_{\infty}^V/F_{\infty}^U)$ because both agree when evaluated at continuous homomorphisms of $U/V$. Hence $(\zeta(F_{\infty}^V/F_{\infty}^U)) \in \Psi_S$. \\
Conversely, assume that $(\zeta(F_{\infty}^V/F_{\infty}^U))$ lie in $\Phi_S$. Take any $f$ in $K_1(\lm(G)_S)$ which maps to $-[C(F_{\infty}/F)]$ under the connecting homomorphism $\partial$. We write $f_{U,V}$ for $\theta_{U,V,S}(f)$. Then $(f_{U,V})$ lies in $\Phi_S$ by C3). On the other hand, under the natural homomorphism
\[
K_0(\lm(G), \lm(G)_S) \ra \prod_{I} K_0(\lm(U/V), \lm(U/V)_S),
\]
$[C(F_{\infty}/F)]$ maps to $([C(F_{\infty}^{V}/F_{\infty}^U)])$. Hence we have 
\[
\partial(f_{U,V}) = -[C(F_{\infty}^V/F_{\infty}^U)],
\]
for every $(U,V) \in I$. From theorem \ref{abelianmainconjecture}, we have that
\[
\partial(\zeta(F_{\infty}^V/F_{\infty}^U)) = -[C(F_{\infty}^V/F_{\infty}^U)],
\]
for every $(U,V)\in I$. Hence $f_{U,V}^{-1}\zeta(F_{\infty}^V/F_{\infty}^U)$ lies in $\lm(U/V)^{\times}$. As $(f_{U,V}^{-1}\zeta(F_{\infty}^V/F_{\infty}^U))$ also lies in $\Phi_S$, C3) gives that 
\[
(f_{U,V}^{-1}\zeta(F_{\infty}^V/F_{\infty}^U)) \in \Phi.
\]
Now C2) gives a $u$ in $K_1(\lm(G))$ such that 
\[
\theta(u) = (f_{U,V}^{-1}\zeta(F_{\infty}^V/F_{\infty}^U)).
\]
Now observe that the natural map from $K_1(\lm(G))$ to $K_1(\lm(G)_S)$ is an injection by our assumption that $\theta$ is an isomorphism. We identify $K_1(\lm(G))$ with its image under this injection. We claim that $uf$ is the $\zeta(F_{\infty}/F)$ which we require. Clearly, $\partial(uf) = -[C(F_{\infty}/F)]$. The interpolation of $L$-values is guaranteed by C1). Let $\rho$ be an Artin character of $G$. We have
\[
\rho = \sum_i n_iInd_{U_i}^G \chi_i.
\]
Then for any positive integer $r$ divisible by $p-1$, we have
\begin{align*} 
\zeta(F_{\infty}/F)(\rho\kappa_F^r) &= \prod_i \zeta(F_{\infty}/F)(Ind_{U_i}^G(\chi_i)\kappa_F^r)^{n_i} \\
                                                             &= \prod_i \zeta(F_{\infty}/F)(Ind_{U_i}^G(\chi_i \kappa_{F_{\infty}^{U_i}}^r))^{n_i} \\
                                                             &= \prod_i (\zeta(F_{\infty}^{V_i}/F_{\infty}^{U_i})(\chi_i \kappa_{F_{\infty}^{U_i}}^r))^{n_i} \\
                                                             &= L_{\Sigma}(\rho, 1-r).
\end{align*}
We now prove the uniqueness of $\zeta(F_{\infty}/F)$. Assume that there is a $\zeta$ which satisfies the conditions of the Main Conjecture. Then $\zeta^{-1} \zeta(F_{\infty}/F)$ lies in $K_1(\lm(G))$. It is easy to check that $\theta_{U,V,S}(\zeta) = \zeta(F_{\infty}^V/F_{\infty}^U)$, for all $(U,V) \in I$, using the uniqueness of the $p$-adic zeta function of Deligne and Ribet. Therefore $\theta(\zeta^{-1}\zeta(F_{\infty}/F)) = (1)$, giving $\zeta(F_{\infty}/F) =\zeta$. \qed \\

\begin{remark} In actual attacks on these problems, $\Phi$, and $\Phi_S$ are described by certain congruences. Hence to prove the noncommutative Main Conjecture we would need to prove certain congruences between abelian $p$-adic zeta functions. This strategy was implemented by Kato in \cite{Kato:2006} to prove the noncommutative Main Conjecture for $p$-adic Heisenberg type extensions. In this paper we predict congruences for extension of the type $G=H\rtimes\Gamma$, with $H$ is a finite abelian $p$ group. However, at present we can prove these congruences only in a special case (see section \ref{sectionspecialcase} for details).
\end{remark}

\subsection{Results proven in this paper}
The strategy proposed in the previous section has two ingredients. Firstly, an algebraic ingredient i.e. existence of some set $I$ such that C1), C2), and C3) are satisfied. Secondly, a number theoretic ingredient i.e. $(\zeta(F_{\infty}^V/F_{\infty}^U))$ lies in $\Phi_S$. In this paper we paper we provide the algebraic ingredient when $G=H\rtimes \G$, with $H$ a finite abelian $p$-group i.e. we show existence of $I$ such that C1), C2), and C3) are satisfied. The subgroups $\Phi$ and $\Phi_S$ are described by certain congruences between the components of elements in $\prod_{I}\lm(U/V)$ and $\prod_{I}\lm(U/V)_S$ respectively. Hence this result predicts congruences between abelian $p$-adic zeta functions $\zeta(F_{\infty}^V/F_{\infty}^U)$ of Deligne and Ribet. We prove these congruences under a technical condition on the group $G$. We call the groups satisfying this condition the groups of \emph{special type} (see definition \ref{specialtype}). Hence we get a proof of the Main Conjecture for extensions of these type. \\

\noindent In section 5 we consider groups $G$ of the form $H\rtimes \G$, with $H$ an abelian compact $p$-adic Lie group which is pro-$p$. We write $G$ as $\varprojlim G/U_k$, with each $U_k$ being an open subgroup of $H$, normal in $G$. Hence each $G/U_k$ is a compact $p$-adic Lie group of dimension 1 and of the form $H/U_k \rtimes \G$. In section 5, using our main theorems from section 3, we prove that validity of the Main Conjecture for $F_{\infty}^{U_k}/F$, for each $k$, implies the Main Conjecture for the extension $F_{\infty}/F$. Combining this result with our proof of the Main Conjecture for extension of special type we get proof of the Main Conjecture in many interesting cases (for example $\mZ_p \rtimes \mZ_p$-extensions). Some examples of groups of special type are discussed in the last section.

\section{Algebraic Preliminaries}
We use the following notations: Recall that $p$ is a fixed odd prime. Let $G$ be a compact $p$-adic Lie group which is pro-$p$. Assume that $G$ contains a finite abelian subgroup $H$ such that $\G := G/H$ is isomorphic to the additive group of $p$-adic integers $\mZ_p$. We write $\G$ multiplicatively and fix a topological generator $\g$ of $\G$. We have an exact sequence 
\[
0 \ra H \ra G \ra \G \ra 0,
\]
which splits. A splitting $\G \ra G$ is given by taking any lifting $\tilde{\g}$ of $\g$ in $G$ and extending it continuously. Also, $\G$ acts on $H$ as $\g \cdot h = \tilde{\g} h \tilde{\g}^{-1}$. As $H$ is abelian this action is independent of the choice of lifting $\tilde{\g}$ of $\g$. This gives an isomorphism of $G$ with $H \rtimes \G$. For any integer $i \geq 0$, we put $\Gi$ for $\G^{p^i}$. As $H$ is a finite group, an open subgroup of $\G$ acts trivially on $H$. We fix such an open subgroup $\G^{(e)}$. For a nonnegative integer $i$, we write $G_i$ for the subgroup $H\rtimes \Gi$ of $G$ and put $H_i$ for $H_0(\Gi, H)=H/(\g^{p^i}-1)H$. Then $G_i^{ab}$, the abelianisation of $G_i$, is given by $H_i \times \Gi$. 

\subsection{Iwasawa algebras and their localisations}
In this section we collect some basic results about Iwasawa algebras of $G_i$. We show that it is enough to localise at a smaller, central, multiplicatively closed subset to get $\lm(G)_S$. We need to take the $p$-adic completion of $\lm(G)_S$, denoted by $\wh{\lm(G)_S}$, in order to define the logarithm. As $H$ is finite, we shall explain that $\lm(G)$ is a \emph{crossed product} of Iwasawa algebra of a central open subgroup and a finite quotient of $G$ (see definition \ref{defcrossedproduct} and lemma \ref{crossedproduct}). We prove a similar result about $\lm(G)_S$ and its $p$-adic completion. In this section $i$ will always denote an integer such that $0 \leq i \leq e$. 

\subsubsection{Ore subsets of Iwasawa algebras}
Recall the subset $S_i$ defined in the introduction:
\[
S_i = \{ s \in \lm(G_i) : \lm(G_i)/\lm(G_i)s \text{  is a finitely generated  } \mZ_p-\text{module} \}.
\]
As mentioned earlier, it is proven in \cite{CFKSV:2005} that $S_i$ is a multiplicatively closed subset of nonzero divisors in $\lm(G_i)$, and is a left and right Ore set. Proposition 2.6 in \emph{loc. cit.} shows that in our case $S_i$ is just the pre-image of all nonzero divisors in $\lm(G_i)/p\lm(G_i)$. Let $\ts$  denote the set $\lm(\G^{(e)}) - p\lm(\G^{(e)})$. $\ts$ is clearly a multiplicatively closed subset of nonzero divisors in $\lm(G_i)$. As $\ts$ is central it is trivially a left and right Ore set. The injection of $\ts$ in $S_i$ gives an injection of $\lm(G_i)_{\ts}$ in $\lm(G_i)_{S_i}$. We show that this map is also surjective.  

\begin{lemma} The natural injection of $\lm(G_i)_{\ts}$ in $\lm(G_i)_{S_i}$ is actually an isomorphism. 
\label{centrals}
\end{lemma}
\noindent{\bf Proof:} Note that $\lm(G_i)_{\ts} = \lm(\Ge)_{\ts} \otimes_{\lm(\Ge)}\lm(G_i)$. We first show that 
\[
Q(\lm(\Ge))\otimes_{\lm(\Ge)}\lm(G_i) \cong Q(\lm(G_i)),
\]
where $Q(R)$ denotes the total ring of fractions of a ring $R$. Note that we have an injective map
\[
Q(\lm(\Ge)) \otimes_{\lm(\Ge)} \lm(G_i) \hookrightarrow Q(\lm(G_i)).
\]
Now $Q(\lm(\Ge))$ is a field and $\lm(G_i)$ is a free $\lm(\Ge)$-module of finite rank (equal to the index of $\Ge$ in $G_i$). Hence $Q(\lm(\Ge)) \otimes_{\lm(\Ge)}\lm(G_i)$ is an Artinian ring. Thus every regular element in it is invertible. $\lm(G_i)$ is contained in $Q(\lm(\Ge))\otimes_{\lm(\Ge)}\lm(G_i)$ and every regular element of $\lm(G_i)$ is invertible in $Q(\lm(\Ge))\otimes_{\lm(\Ge)}\lm(G_i)$. Hence $Q(\lm(\Ge))\otimes_{\lm(\Ge)}\lm(G_i) \hookrightarrow Q(\lm(G_i))$ must be surjective. Hence any element $x$ in $\lm(G_i)_{S_i} \subset Q(\lm(G_i))$ can be written as $\frac{a}{t}$ with $a$ in $\lm(G_i)$ and $t$ a nonzero element in $\lm(\Ge)$. If $t$ lies in $p^n\lm(G_i)$, then $tx = a$ lies in $p^n\lm(G_i)_{S_i}$. On the other hand $a$ also lies in $\lm(G_i)$. Hence $a$ lies in $p^n\lm(G_i)$. Thus we can divide the largest possible power of $p$ from $t$ and the same power of $p$ from $a$ and $x$ can be represented as $\frac{a}{t}$, with $a$ in $\lm(G_i)$ and $t$ in $\ts$. \qed \\
 
\noindent We recall the $p$-adic completion of $\lm(G)_{\ts}$,
\[
\wh{\lm(G)_{\ts}} := \ilim{n} \lm(G)_{\ts}/p^n\lm(G)_{\ts}.
\]
We denote the Jacobson radical of any ring $R$ by $J(R)$. The following result in proven in \cite{CFKSV:2005}. 

\begin{corollary} Some power of $J(\lm(G)_{\ts})$ (resp. $J(\wh{\lm(G)_{\ts}})$) is contained in $p\lm(G)_{\ts}$ (resp. $p\wh{\lm(G)_{\ts}}$). 
\label{jisnilpotent} 
\end{corollary}
\noindent{\bf Proof:} We prove both assertion simultaneously. For a pro-finite group $P$, let $\mathbb{F}_p[[P]]$ denote the Iwasawa algebra of $P$ with coefficients in $\mathbb{F}_p$, i.e. $\mathbb{F}_p[[P]] = \varprojlim \mathbb{F}_p[P/U]$, where $U$ runs through the set of open normal subgroup of $P$. We have the following short exact sequence
\[
0 \ra \frac{J(\lm(G)_{\ts})}{p\lm(G)_{\ts}} = \frac{J(\wh{\lm(G)_{\ts}})}{p\wh{\lm(G)_{\ts}}} \ra Q(\mathbb{F}_p[[G]]) \ra Q(\mathbb{F}_p[[\G]]) \ra 0.
\]
Write $J$ for $\frac{J(\lm(G)_{\ts})}{p\lm(G)_{\ts}} = \frac{J(\wh{\lm(G)_{\ts}})}{p\wh{\lm(G)_{\ts}}}$. We must prove that $J$ is a nilpotent ideal. Let $N$ be kernel of the natural map from $\mathbb{F}_p[[G]]$ to $\mathbb{F}_p[[\G]]$. Let $I_H$ be the kernel of the augmentation map from $\mathbb{F}_p[H]$ to $\mathbb{F}_p$. Then $N$ is $\mathbb{F}_p[[G]]I_H$ and as $H$ is a finite $p$-group,
\[
N^n = \mathbb{F}_p[[G]]I_H^n = 0,
\]
for some positive integer $n$. By the lemma, we can write any element $x$ of $Q(\mathbb{F}_p[[G]])$ as $\frac{a}{t}$ with $a$ in $\mathbb{F}_p[[G]]$ and $t$ is a nonzero element in $\mathbb{F}_p[[\G]]$. It is clear that $x$ lies in $J$ if and only if $a$ lies in $N$. Since $t$ is central, we deduce that $J$ is nilpotent. \qed \\

\subsubsection{Iwasawa algebras as crossed products}
\begin{definition} Let $R$ be a ring and $\Pi$ be a finite group. Then a \emph{crossed product} $R\star \Pi$ of $\Pi$ over $R$ is an associative ring which contains $R$ and has as an $R$ basis the set $\bar{\Pi}$, a copy of $\Pi$. Thus each element of $R\star \Pi$ is uniquely expressed as a finite sum $\sum_{\pi \in \Pi}r_{\pi}\bar{\pi}$ with $r_{\pi} \in R$. Addition is defined in the obvious way and multiplication is defined by following rule: for $\pi_1$ and $\pi_2$ in $\Pi$ we have
\[
\bar{\pi}_1\bar{\pi}_2=\tau(\pi_1,\pi_2)\overline{\pi_1\pi_2}
\]
where $\tau:\Pi\times\Pi \ra R^{\times}$ is called \emph{twisting map}. We assume that $R$ commutes with elements of $\bar{\Pi}$. 
\label{defcrossedproduct}
\end{definition}   

\noindent We now express $\lm(G)$ (resp. $\lm(G)_{\ts}$ and $\wh{\lm(G)_{\ts}}$) as crossed product of $G/\Ge$ over the ring $\lm(\Ge)$ (resp. $\lm(\Ge)_{\ts}$ and $\wh{\lm(\Ge)_{\ts}}$). To this end we fix the following set of coset representatives of $\Ge$ in $G$
\[
C=\{ h\g^a : \text{  with  } h\in H \text{  and  } 0 \leq a \leq p^e-1 \}.
\]
There is a bijection between $C$ and $G/\Ge$. $C$ is our copy $\overline{G/\Ge}$ of $G/\Ge$ in the definition of crossed product. Now the twisting map,
\[
\tau : G/\Ge \times G/\Ge \ra \Ge \subset \lm(\Ge)^{\times} \subset \lm(\Ge)_{\ts}^{\times} \subset \wh{\lm(\Ge)_{\ts}}^{\times},
\]
is obvious. $(h_1\g^{a_1}, h_2\g^{a_2}) \mapsto \g^{a_1+a_2 - [a_1+a_2]} \in \Ge$. Here $[k]$ is the smallest nonnegative integer such that $k \equiv [k] (mod \ p^e)$. Note that $\tau(\pi_1, \pi_2) = \tau(\pi_2, \pi_1)$ i.e. $\tau$ is \emph{symmetric}.  

\begin{lemma} With the twisting map as above, $\lm(G)$ (resp. $\lm(G)_{\ts}$ and $\wh{\lm(G)_{\ts}}$) is isomorphic to the crossed product $\lm(\Ge)\star(G/\Ge)$ (resp. $\lm(\Ge)_{\ts}\star(G/\Ge)$ and $\wh{\lm(\Ge)_{\ts}}\star (G/\Ge)$).
\label{crossedproduct}
\end{lemma}
\noindent{\bf Proof:} Note that $\mZ_p[G/\G^{(k)}]$ is isomorphic to the crossed product $\mZ_p[\Ge/\G^{(k)}]\star (G/\Ge)$, for every integer $k \geq e$ and these isomorphisms are compatible with natural projections. Hence the lemma. \qed \\

\begin{remark} Similar result holds for $G_i$ and $G_i^{ab}$. Analogous statement for these groups is clear.
\end{remark}

\noindent We introduce another object that we require in section \ref{sectionlogarithm}. For any ring $R$, let $[R,R]$ denote the additive subgroup of $R$ generated by elements of the form $ab-ba$, for all $a$ and $b$ in $R$. We write $T(R)$ for $R/[R,R]$. For a finite group $P$, we put $Conj(P)$ for the set of conjugacy classes of $P$. We need the following trivial 

\begin{lemma} $T(\lm(G))$ (resp. $T(\lm(G)_{\ts})$ and $T(\wh{\lm(G)_{\ts}})$) is a free $\lm(\Ge)$ (resp. $\lm(\Ge)_{\ts}$ and $\wh{\lm(\Ge)_{\ts}}$) module of finite rank with a basis given by the set of conjugacy classes of $G/\Ge$.
\label{basisoft}
\end{lemma}
\noindent{\bf Proof:} Note that $a \equiv bab^{-1} (mod \ [\lm(G), \lm(G)])$ for any $a$ and $b$ in $G$. Take any $x$ in $\lm(G)$. Write it as $\sum_{h \in G/\Ge} x_hh$. Therefore image of $x$ in $T(\lm(G))$, denoted by $\bar{x}$ is $\sum_{C\in Conj(G/\Ge)} (\sum_{h\in C}x_h)C$. Hence conjugacy classes of $G/\Ge$ generates $T(\lm(G))$. Now we prove the linear independence of conjugacy classes. If $\sum_{C\in Conj(G/\Ge)} y_CC=0$ in $T(\lm(G))$, then 
\[
\sum_{C\in Conj(G/\Ge)}y_Ch_C = \sum z_i(a_ib_i-b_ia_i)
\]
in $\lm(G)$, where $h_C$ is any element in $C$. Now notice that we may assume that $a_i$ and $b_i$ are in $G/\Ge$ because $\tau$ is symmetric. Then it is easily seen that $y_C=0$ for every $C$. The proof for $T(\lm(G)_{\ts})$ and $T(\wh{\lm(G)_{\ts}})$ is exactly the same. 
\qed \\

\subsubsection{Relations in subquotients of $G$}Results in this subsection are essentially trivial but we collect them for convenience of use later.

\begin{lemma} Two elements $h_1\g^{a_1}$ and $h_2\g^{a_2}$ of $G$ are conjugates if and only if the following holds: \\
1)  $a_1=a_2=:a$ and \\
2) if $a=0$ then there is a $k\in\mZ_p$ such that $h_2=\g^kh_1\g^{-k}$ or if $a \in p^i\mZ_p$ but $a\notin p^{i+1}\mZ_p$, then there is a $k\in \{0, 1, 2, \ldots, (p^i-1)\}$ such that $h_2\g^kh_1^{-1}\g^{-k} \in (\g^{p^i}-1)H$.
\label{conjugatesofg}
\end{lemma}
\noindent{\bf Proof:} Any element of $G$ can be written as $h\g^{b}$ with $h\in H$ and $b\in \mZ_p$. Then
\[
h\g^{b}(h_1\g^{a_1})\g^{-b}h^{-1}=(\g^bh_1\g^{-b})(h\g^{a_1}h^{-1}).
\]
Let $x=h\g^{a_1}h^{-1}$. Then 
\[
x\g^{-a_1}=h(\g^{a_1}h^{-1}\g^{-a_1})=\g^{a_1}h^{-1}\g^{-a_1}h.
\]
Hence $x \in (\g^{p^i}-1)H\g^{a_1}$. For the first term, write $b=k+p^ib_1$, with $k\in \{0,1,\ldots, p^i-1\}$ and $b_1 \in \mZ_p^{\times}$. Then 
\[
\g^bh_1\g^{-b} = \g^{b_1p^i}(\g^kh_1\g^{-k})\g^{-b_1p^i}=(\g^kh_1\g^{-k})[\g^{b_1p^i}(\g^kh_1\g^{-k})\g^{-b_1p^i}(\g^kh_1^{-1}\g^{-k})] \in (\g^kh_1\g^{-k})(\g^{p^i}-1)H.
\]
Hence the lemma. \qed 

\begin{lemma} Let $a \in \mZ_p$ be such that $a \in p^i\mZ_p-p^{i+1}\mZ_p$. Then for any $h \in H$, 
\[
[h(\g^ah\g^{-a})(\g^{2a}h\g^{-2a})\cdots(\g^{(p-1)a}h\g^{-(p-1)a})] = 
[h(\g^{p^i}h\g^{-p^i})(\g^{2p^i}h\g^{-2p^i})\cdots(\g^{(p-1)p^i}h\g^{-(p-1)p^i})]
\]
in $G_{i+1}^{ab}$. Here and hereafter $[\cdot]$ will denote the equivalence class in the relevant space. 
\label{relationing}
\end{lemma}
\noindent{\bf Proof:} $a$ and $p^i$ generate the subgroup of order $p$ in $\mZ/p^{i+1}\mZ$. Hence for every $0 \leq i \leq p-1$, there is a unique $0\leq j_i \leq p-1$ such that $ia \equiv j_ip^i (mod \ p^{i+1})$. The result follows as $G_{i+1}^{ab}=(H/(\g^{p^{i+1}}-1)H)\times\G^{(i+1)}$. \qed \\

\begin{lemma} Let $\rho$ be any irreducible Artin representation of $G$ (i.e. which factors through a finite quotient of $G$). Then $\rho$ is obtained by inducing a one dimensional representation of $H \rtimes \Gi$, for some $0\leq i \leq e$. 
\label{representationsofg}
\end{lemma}
\noindent{\bf Proof:} Assume that $\rho$ factors through $G/ \G^{(f)} \cong H\rtimes \G/\G^{(f)}$. We may and do assume that $f \geq e$. Then it is proven in Serre \cite{Serre:representationtheory}, proposition 25, that the representation $\rho$ of $G/\G^{(f)}$ is induced from a one dimensional representation of $G \rtimes \G^{(i)}/\G^{(f)}$ for some $0 \leq i \leq e$. \qed \\

\subsection{Logarithm and integral logarithm on $K_1$} \label{sectionlogarithm}
The first aim of this section is to construct the logarithm and the integral logarithm maps on $K_1$ groups of $\lm(G)$ and $\wh{\lm(G)_{\ts}}$. This is a natural extension of the integral logarithm of Oliver and Taylor to the setting of Iwasawa theory. The first use of integral logarithm in non-commutative Iwasawa theory is due to Ritter and Weiss \cite{RitterWeiss:3}. Of course, the integral logarithm was first used in the commutative Iwasawa theory by Coates-Wiles for defining the ``Coates-Wiles homomorphism" (see Washington \cite{Washington:cyclotomicfields}). The use of integral logarithm was suggested to us by Kato. In the second part of this section we study the behaviour of integral logarithm with respect to norm homomorphism on $K_1$ groups. We follow Oliver and Taylor \cite{OliverTaylor:1988} for this. We generalise theorem 1.4 in \cite{OliverTaylor:1988} to our setting. In this section we crucially use that $\tau$ is symmetric. 

\subsubsection{Logarithm on $K_1(\lm(G))$ and $K_1(\wh{\lm(G)_{\ts}})$}
For this subsection $R$ will be either $\lm(\Ge)$ or $\wh{\lm(\Ge)_{\ts}}$. In this section, following R. Oliver (\cite{Oliver:1988}, chapter 2), we construct logarithm on $K_1(R\star(G/\Ge))$. We introduce some notations. Throughout this subsection let $\R$ denote the ring $R \star (G/\Ge)$ and let $J$ be the Jacobson radical of $\R$. We put $\R_{\mQ}$ for $\mQ\otimes \R$ and put $I_{\mQ}$ for $\mQ\otimes I$, for any ideal $I$ of $\R$. Let $[\R_{\mQ}, I_{\mQ}]$ be the additive subgroup generated by elements of the form $ab-ba$ for all $a$ in $\R_{\mQ}$ and all $b$ in $I_{\mQ}$. Recall $GL_{\infty}(\R)$ is the union $\cup_{n=1}^{\infty}GL_n(\R)$. Let $E(\R, I)$ be the smallest normal subgroup of $GL_{\infty}(\R)$ generated by elementary matrices $E_{i,j}^r$ (i.e $r$ in the $i,j$th position, with $i \neq j$, all 1's on diagonal and 0 everywhere else), with $ r \in I$. We write $GL_{\infty}(\R, I)$ for kernel of the natural map from $GL_{\infty}(\R)$ to $GL_{\infty}(\R/I)$. We define $K_1(\R, I)$ to be the quotient $GL_{\infty}(\R,I)/E(\R,I)$. It is a well known fact from algebraic $K$-theory that this is an abelian group (see, for example, \cite{Oliver:1988}). Let
\[
Log(1+X) = \sum_{n=1}^{\infty} (-1)^{n-1} \frac{X^n}{n}, \text{    and}
\]
\[
Exp(X) = \sum_{n=0}^{\infty} \frac{X^n}{n!}.
\]

\begin{lemma} (1) Let $I \subset J$ be an ideal of $\R$. Then for any $u, v$ in $1+I$, $Log(u)$ and $Log(v)$ converge in $I_{\mQ}$ and 
\[
Log(uv) \equiv Log(u) + Log(v) (mod \ [\R_{\mQ}, I_{\mQ}]).
\]
(2) Assume that $I$ is an ideal of $\R$ contained in $p\R$. Then for any $u,v$ in $1+I$, $Log(u)$ and $Log(v)$ converge in $I$ and, 
\[
Log(uv) \equiv Log(u) + Log(v) (mod \ [\R, I]).
\]
(3) Assume again that $I$ is an ideal of $\R$ contained in $p\R$. Then for any $x, y$ in $I$, $Exp(x)$ converges in $1+I$. $Log$ and $Exp$ are inverse bijections between $I$ and $1+I$. Moreover, $Exp([\R, I]) \subset E(\R, I)$ and 
\[
Exp(x+y) \equiv Exp(x)Exp(y) (mod \ E(\R,I)).
\]
\label{technicallog}
\end{lemma}
\noindent{\bf Sketch Proof:} (For details see lemma 2.7 in \cite{Oliver:1988}) $R$ is a local ring. Put $\mathfrak{m}$ for the maximal ideal of $R$. $R$ is complete with respect to $\mathfrak{m}$-adic topology. Then $J/\mathfrak{m}^n\R$ is a nilpotent ideal in $\R/\mathfrak{m}^n\R$. Therefore, $u^n/n$ converges to 0 as $n$ tends to infinity. Hence $Log(u)$ and $Log(v)$ converge in $I_{\mQ}$. When $ I \subset p\R$, it is easy to see that $I^n \subset n!I^k$, with $k$ going to infinity with $n$. Hence $Log(u)$ converges in $I$ for any $u \in 1+I$ and $Exp(x)$ converges in $1+I$ for any $x$ in $I$. We now show the congruences. We set
\[
U(I) = \sum_{m,n \geq 1}\frac{1}{m+n}[I^m,I^n] \subset [\R_{\mQ}, I_{\mQ}].
\]
we will show that 
\[
Log((1+x)(1+y)) \equiv Log(1+x) + Log(1+y) \ (mod \ U(I)),
\]
for any $x,y$ in $I$. When $I \subset p\R$, $U(I) \subset [\R,I]$, hence this will show both (1) and (2). For each $n \geq 1$, we let $W_n$ be the set of formal ordered monomials of length $n$ in two variables $a,b$. For $w \in W_n$, set 

$C(w)$ = orbit of $w$ in $W_n$ under cuclic permutations. 

$k(w)$ = number of occurrences of $ab$ in $w$.

$r(w)$ = coefficients of $w$ in $Log(1+a+b+ab) = \sum_{i=0}^{k(w)} (-1)^{n-i-1} \frac{1}{n-i}\binom{k(w)}{i}$.

If $w^{\prime} \in C(w)$, then it is clear that $w(x,y) \equiv w^{\prime}(x,y) (mod \ [I^i, I^j])$ for some $i,j \geq 1$ such that $i+j=n$. So
\[
Log(1+x+y+xy) = \sum_{n=1}^{\infty} \sum_{w\in W_n} r(w)w(x,y) 
\equiv \sum_{n=1}^{\infty} \sum_{w \in W_n/C} \Big( \sum_{w^{\prime} \in Cw} r(w^{\prime})\Big)w(x,y) (mod \ U(I)).
\]
Let $ k =max \{k(w^{\prime}) : w^{\prime} \in C(w)\}$. Let $ |C(w)| = n/t$. Then $C(w)$ contains $k/t$ element with exactly $(k-1)$ $ab$'s and $(n-k)/t$ elements with $k$ $ab$'s. Hence
\begin{align*}
\sum_{w^{\prime} \in C(w)} r(w^{\prime}) &= \frac{1}{t} \sum_{i=0}^{k} (-1)^{n-i-1} \frac{1}{n-i}\Big((n-k)\binom{k}{i} + k \binom{k-1}{i}\Big) \\
                                                                        &=\frac{1}{t} \sum_{i=0}^{k} (-1)^{n-i-1} \frac{1}{n-i}\Big((n-k)\binom{k}{i} + (k-i)\binom{k}{i}\Big) \\
                                                                        &= \frac{1}{t}\sum_{i=0}^{k} (-1)^{n-i-1} \binom{k}{i},
\end{align*}
which is 0 unless $k=0$, in which case it is equal to $(-1)^{n-1}\frac{1}{n}$. Thus
\[
Log(1+x+y+xy) \equiv \sum_{n=1}^{\infty} (-1)^{n-1} \Big(\frac{x^n}{n}+\frac{y^n}{n}\Big) = Log(1+x)+Log(1+y) \ (mod\ U(I)).
\]
 
\noindent We now prove the remaining part of (3). $Exp$ and $Log$ induce bijection between $I$ and $1+I$. Hence 
\[
Log(Exp(x)Exp(y)) \equiv x+y \ (mod \ U(I)),
\]
which gives $Exp(x)Exp(y)Exp(x+y)^{-1} \subset Exp(U(I)) \subset Exp([\R,I])$. Hence we only need to prove that $Exp([\R,I])$ is contained in $E(\R,I)$. Choose a $R$-basis $\{[s_1,v_1], \ldots, [s_m,v_m]\}$ of $[\R, I]$, with $s_i \in \R$ and $v_i \in I$. Let $x = \sum_{i=1}^{m}a_i[s_i,v_i]$ be an element in $[\R,I]$. Define 
\[
\psi(x)= \prod_{i=1}^{m}(Exp(a_is_iv_i)Exp(a_iv_is_i)^{-1}).
\]
By Vaserstein's identity (\cite{Oliver:1988}, theorem 1.15), $\psi(x)$ lies in $E(\R,I)$. For any $u$ in $Exp(p[\R,I])$, we define a sequence $x_0, x_1, \ldots, $ given by 
\[
x_0 = Log(u) \in p[\R,I]; \hspace{2cm} x_{i+1} = x_i +Log(\psi(x_i)^{-1}u).
\]
Then 
\[
\psi(x_i) \equiv u , \hspace{2cm} x_{i+1}\equiv x_i (mod \ p^{2+i}[\R,I]).
\]
Hence $u = \psi(lim_{i\ra \infty}x_i) \in E(\R,I)$. And it can be shown that $Exp([\R,I]) \subset E(\R,I)Exp(p[\R,I])$. (See \cite{Oliver:1988}, lemma 2.7 for details). \qed \\ 

\begin{theorem} Let $I \subset J$ be any ideal of $\R$. Then $Log(1+x)$ induces a unique homomorphism 
\[
log : K_1(\R, I) \ra \mQ \otimes_{\mZ} (I/[\R,I]).
\]
Furthermore, if $I \subset p\R$, then the logarithm induces an isomorphism 
\[
log : K_1(\R, I) \ra I/[\R,I].
\]
\label{constructionoflog}
\end{theorem}
\noindent{\bf Proof:} This is a formal consequence of the above lemma (See \cite{Oliver:1988} theorem 2.8). \qed \\

\begin{remark} we have constructed $log$ homomorphism on the groups $K_1(\lm(G), J(\lm(G)))$ and $K_1(\wh{\lm(G)_{\ts}}, J(\wh{\lm(G)_{\ts}}))$, but its clear that $G$ can be replaced by any group $G_i$ or $G_i^{ab}$ or $\Gi$. In particular, we have an isomorphism 
\[
log:1+p\wh{\lm(\G)_{\ts}} \ra p\wh{\lm(\G)_{\ts}}.
\]
\end{remark}

\begin{remark} Since $\lm(G)$ and $\wh{\lm(G)_{\ts}}$ are local rings, the following natural maps are all surjective
\begin{center}
$\lm(G)^{\times} \ra K_1(\lm(G)),$ \\
$\wh{\lm(G)_{\ts}^{\times}} \ra K_1(\wh{\lm(G)_{\ts}}),$ \\
$1+J(\lm(G)) \ra K_1(\lm(G), J(\lm(G)))$ \\
$1+J(\wh{\lm(G)_{\ts}}) \ra K_1(\wh{\lm(G)_{\ts}}, J(\wh{\lm(G)_{\ts}})).$
\end{center}
\noindent We have defined $log$ on relative $K_1$'s and they fit into the following exact sequences:
\begin{center}
$K_2(\lm(G)/J(\lm(G))) \ra K_1(\lm(G), J(\lm(G))) \ra K_1(\lm(G)) \ra K_1(\lm(G)/J(\lm(G))),$ and \\
$K_2(\wh{\lm(G)_{\ts}}/J(\wh{\lm(G)_{\ts}})) \ra K_1(\wh{\lm(G)_{\ts}}, J(\wh{\lm(G)_{\ts}})) \ra K_1(\wh{\lm(G)_{\ts}}) \ra K_1(\wh{\lm(G)_{\ts}}/J(\wh{\lm(G)_{\ts}})).$
\end{center}
\noindent Now $\lm(G)/J(\lm(G)) \cong \mathbb{F}_p$; it is well known that $K_2(\mathbb{F}_p)$ is trivial and $K_1(\mathbb{F}_p) \cong \mathbb{F}_p^{\times}$. Hence $log$ can be uniquely extended to $K_1(\lm(G))$. On the other hand, $\wh{\lm(G)_{\ts}}/J(\wh{\lm(G)_{\ts}}) \cong Q(\mathbb{F}_p[[\G]])$; and $K_1(Q(\mathbb{F}_p[[\G]])) \cong Q(\mathbb{F}_p[[\G]])^{\times}$ which is not a torsion group. So it is not possible to extend $log$ to $K_1(\wh{\lm(G)_{\ts}})$. However, it is known that $K_2(Q(\mathbb{F}_p[[\G]]))$ is torsion. Hence its image in $K_1(\wh{\lm(G)_{\ts}}, J(\wh{\lm(G)_{\ts}}))$ lies in the kernel of $log$. This gives an extension of $log$ to elements of $K_1(\wh{\lm(G)_{\ts}})$ which lie in the image of $1+J(\wh{\lm(G)_{\ts}})$ under the natural map
\[
1+J(\wh{\lm(G)_{\ts}}) \hookrightarrow \wh{\lm(G)_{\ts}^{\times}} \ra K_1(\wh{\lm(G)_{\ts}}).
\]
\end{remark}

\subsubsection{Integral logarithm} 
\begin{definition} For any group $P$, we define $\varphi$ to be the map which takes $g$ in $P$ to $g^p$. 
\end{definition}
\noindent Note that this map need not always be a homomorphism. However, if $P$ is abelian then it is a homomorphism. It is easy to see that $\varphi$ on the group $G$ takes the subset $\ts$ to itself. Hence it induces a map, which we again denote by $\varphi$, on $\lm(G)_{\ts}$ and $\wh{\lm(G)_{\ts}}$. If $G$ is abelian, $\varphi$ is a ring homomorphism. 

\begin{lemma} For any $u$ in $\wh{\lm(\G)_{\ts}^{\times}}$, $\frac{u^p}{\vp(u)}$ lies in $1+p\wh{\lm(\G)_{\ts}}$.
\end{lemma}
\noindent{\bf Proof:} Mapping $\g$ to $X+1$ gives an isomorphism between $\lm(\G)$ and $\mZ_p[[X]]$, the ring of power series in $X$ with coefficient in $\mZ_p$. This isomorphism extends to an isomorphism between $\lm(\G)_{\ts}$ and $\mZ_p[[X]]_{(p)}$, the localisation of $\mZ_p[[X]]$ at the prime ideal $(p)$ and also to an isomorphism between $\wh{\lm(\G)_{\ts}}$ and $\wh{\mZ_p[[X]]_{(p)}}$, the $p$-adic completion of $\mZ_p[[X]]_{(p)}$. $\vp$ induces a map on $\mZ_p[[X]]_{(p)}$ which maps $f(X)$ to $f((1+X)^p-1)$. Then it is clear that $\vp(u) \cong u^p \ (mod \ p)$ for any $u$ in $\lm(\G)_{\ts}^{\times}$. As $\vp$ is a continuous map the same is deduced for any $u$ in $\wh{\lm(\G)_{\ts}^{\times}}$. \qed \\

\begin{definition} Let $x$ be an element of $\lm(G)^{\times}$. Define
\begin{center}
$L: K_1(\lm(G)) \ra \mQ_p\otimes_{\mZ_p}T(\lm(G)),$
\end{center}
given by $[x] \mapsto log([x]) - \frac{\vp}{p}log([x])$. Here $[x]$ denotes the class of $x$ in $K_1(\lm(G))$. 
\end{definition}

\begin{definition} Let $x$ be any element in $\wh{\lm(G)_{\ts}}^{\times}$. Then it can be written as $u(1+y)$ with $u$ in $\wh{\lm(\G)_{\ts}}^{\times}$ and $y$ in $J(\wh{\lm(G)_{\ts}})$. Define
\begin{center}
$L: K_1(\wh{\lm(G)_{\ts}}) \ra \mQ_p\otimes_{\mZ_p}T(\wh{\lm(G)_{\ts}}),$
\end{center}
given by $[x] \mapsto \frac{1}{p}log(\frac{u^p}{\vp(u)}) + log([1+y]) - \frac{\vp}{p}log([1+y])$. 
\end{definition}

\noindent It is easily seen that $L([x])$ is independent of the choices of $u$ and $y$. The following lemma is a slight generalisation of the result of Oliver \cite{Oliver:1988}. We sketch its proof here. For details see chapter 6 in \emph{loc. cit.}.

\begin{lemma} $L$ is a homomorphism from $K_1(\lm(G))$ to $T(\lm(G))$ (resp. $K_1(\wh{\lm(G)_{\ts}})$ to $T(\wh{\lm(G)_{\ts}})$). Furthermore, $L$ is natural with respect to ring homomorphisms induced by group homomorphisms.
\label{constructionofintegrallog}
\end{lemma}
\noindent{\bf Proof:} Till the end of this proof let $R$ denote either $\lm(\Ge)$ or $\wh{\lm(\Ge)_{\ts}}$ and let $\R$ denote either $\lm(G)$ or $\wh{\lm(G)_{\ts}}$ and let $J$ denote the Jacobson radical of $\R$. Let $x$ be any element in $J$. Then
\begin{align*}
L(1-x) &= -\big[ x+\frac{x^2}{2} + \frac{x^3}{3} +\cdots\big] + \big[ \frac{\varphi(x)}{p}+\frac{\varphi(x^2)}{2p}+\cdots \big] \\
            &\equiv -\sum_{k=1}^{\infty} \frac{1}{pk}[x^{pk}-\varphi(x^k)] \ \ (mod \ T(\R)).
\end{align*}
Hence we must prove that $p^n | [x^{p^n}-\varphi(x^{p^{n-1}})]$. Write $x= \sum r_ig_i$ with $r_i \in R$ and $g_i \in G/\Gamma^{(e)}$. Set $q=p^n$. A typical term of $x^q$ looks like $r_{i_1}\cdots r_{i_q} g_{i_1} \cdots g_{i_q}$. $\mZ/p^n$ cyclically permutes $g_i$'s. If there are $p^t$ cyclic permutations which leave it invariant, then there are total of $p^{n-t}$ conjugates appearing in $x^q$. Then $g_{i_1}\cdots g_{i_q}$ is a $p^t$th power and the sum of conjugates has the form $p^{n-t}\hat{r}^{p^t}\hat{g}^{p^t}$ in $T(\R)$. If $t=0$ then this is a multiple of $p^n$. If $t \geq 1$, then there is a corresponding term $p^{n-t}\hat{r}^{p^{t-1}}\hat{g}^{p^{t-1}}$ in $x^{p^{n-1}}$. So
\[
p^{n-t}\hat{r}^{p^t}\hat{g}^{p^t} \equiv p^{n-t}\varphi(\hat{r}^{p^{t-1}})\hat{g}^{p^t} \equiv p^{n-t}\varphi(\hat{r}^{p^{t-1}}\hat{g}^{p^{t-1}}) \ \ (mod \ p^n).
\]
Here, the first congruence comes from the previous lemma. In this proof we crucially use the fact that $\tau$ is symmetric. \qed \\

\noindent The following lemma for $\lm(G)$ can be very easily deduced from theorem 6.6 in Oliver \cite{Oliver:1988} for example, by passing to the inverse limit of $\mZ_p[G/U]$ over all open normal subgroups $U$ of $G$. 

\begin{lemma} We have an exact sequence 
\[
1 \ra G^{ab} \times \mu_{p-1} \ra K_1(\lm(G)) \xrightarrow{L} T(\lm(G)) \ra G^{ab} \ra 1.
\]
Here, $\mu_{p-1}$ denotes the group of $p-1$st roots of unity. 
\label{exactsequenceofintegrallog}
\end{lemma}
\noindent{\bf Proof:} See theorem 6.6 in \cite{Oliver:1988} 

\subsubsection{Integral logarithm and norm homomorphism}
Recall that we have reserved the notation $i$ for an integer such that $0 \leq i \leq e$ and that $[\cdot]$ denotes an equivalence class in the relevant space. Since $\lm(G)$ (resp. $\wh{\lm(G)_{\ts}}$) is a free $\lm(G_i)$ (resp. $\wh{\lm(G_i)_{\ts}}$) module of rank $p^i$, we have norm homomorphisms:
\[
N_i: K_1(\lm(G)) \ra K_1(\lm(G_i))
\]
\[
(\text{resp.  } N_i: K_1(\wh{\lm(G)_{\ts}}) \ra K_1(\wh{\lm(G_i)_{\ts}})).
\]
On the other hand, on additive side, we have homomorphisms:
\[
res_i: T(\lm(G)) \ra T(\lm(G_i))
\]
\[
(\text{resp.  } res_i : T(\wh{\lm(G)_{\ts}}) \ra T(\wh{\lm(G_i)_{\ts}})),
\]
defined as follows; using lemma \ref{basisoft}, we define $res_i$ on $\lm(\Ge)$ basis of $T(\lm(G))$ (resp. $\wh{\lm(\Ge)_{\ts}}$ basis of $T(\wh{\lm(G)_{\ts}})$) and extend linearly. Define
\[
res_i([h\g^a]) = 
\begin{cases}
\sum_{k=0}^{p^i-1}[\gamma^k(h\gamma^a)\gamma^{-k}] \ \ \ \ \ \ \ \ \ \ if \ p^i|a. \\
0 \hspace{3.3cm} otherwise.
\end{cases}
\]

\begin{proposition} For any $y$ in $J(\lm(G))$ (resp. $y$ in $J(\wh{\lm(G)_{\ts}})$), we have:
\[
res_i(log([1+y])) = log(N_i([1+y])),
\]
in both cases.
\label{normandres}
\end{proposition}
\noindent{\bf Proof:} The proof of theorem 1.4 in \cite{OliverTaylor:1988} goes through as it is given there. Here we use the crossed product description (lemma \ref{crossedproduct}) and also the fact that $J(\widehat{\lm(G)_{\ts}})$ is nilpotent modulo $p$ (corollary \ref{jisnilpotent}). \qed \\

\noindent Recall the map $\vp$ on $\lm(G)$ and on $\wh{\lm(G)_{\ts}}$ induced by the the map $g \mapsto g^p$ on $G$. It clearly induces a map on $T(\lm(G))$ and on $T(\wh{\lm(G)_{\ts}})$, which we again denote by $\vp$. We define homomorphisms
\[
\beta_i : T(\lm(G)) \ra T(\lm(G_i)) \ra \lm(G_i^{ab}), \text{    and}
\]
\[
\hat{\beta}_i : T(\wh{\lm(G)_{\ts}}) \ra T(\wh{\lm(G_i)_{\ts}}) \ra \wh{\lm(G_i^{ab})_{\ts}},
\]
where, in both cases, the first homomorphism in $res_i$ and second is the one induced by natural surjection of $G_i$ on $G_i^{ab}$. We also define:
\[
\theta_i : K_1(\lm(G)) \ra K_1(\lm(G_i)) \ra \lm(G_i^{ab})^{\times}, \text{    and}
\]
\[
\hat{\theta}_i : K_1(\wh{\lm(G)_{\ts}}) \ra K_1(\wh{\lm(G_i)_{\ts}}) \ra \wh{\lm(G_i^{ab})_{\ts}}^{\times},
\]
where, in both cases, the first map in $N_i$ and the second is the one induced by natural surjection of $G_i$ on $G_i^{ab}$. \\

\noindent For any $p$-adic character $\chi$ of $G_i^{ab}$, we define $\tilde{\chi}$ to the the automorphism of $\bar{\mZ}_p\otimes_{\mZ_p}\lm(G_i^{ab})$ or of $\bar{\mZ}_p\otimes_{\mZ_p}\wh{\lm(G_i^{ab})_{\ts}}$ induced by $g \mapsto \chi(g)g$. For $i \geq 1$ we define $ver_i$ to be the homomorphism from $\lm(G_{i-1}^{ab})$ to $\lm(G_i^{ab})$ (or from $\wh{\lm(G_{i-1}^{ab})_{\ts}}$ to $\wh{\lm(G_i^{ab})_{\ts}}$) induced by the transfer homomorphism from $G_{i-1}^{ab}$ to $G_i^{ab}$. For every $i \geq 1$, we fix a $\omega_i$ to be a non-trivial $p$-adic character of $\G^{(i-1)}$ such that $\omega_i|_{\Gi}$ is trivial. We will consider it as a character of $G_{i-1}^{ab}$ using the natural surjection of $G_{i-1}^{ab}$ onto $\G^{(i-1)}$.

\begin{lemma} For any $i \geq 1$ and any $a \in T(\lm(G))$ (resp. $a \in T(\wh{\lm(G)_{\ts}})$)
\[
\beta_i\vp(a) - \vp\beta_i(a) = ver_i\Big(p\beta_{i-1}(a)-\sum_{k=0}^{p-1}\tilde{\omega}_i^k\beta_{i-1}(a)\Big).
\]
\[
(\text{resp.} \ \ \ \   \hat{\beta}_i\vp(a) - \vp\hat{\beta}_i(a) = ver_i\Big(p\hat{\beta}_{i-1}(a)-\sum_{k=0}^{p-1}\tilde{\omega}_i^k\hat{\beta}_{i-1}(a)\Big)).
\]
\label{beta1}
\end{lemma}    
\noindent{\bf Proof:} Let $a = \sum a_hh$. Then we have
\[
\beta_i\varphi(a) = \sum \varphi(a_h) \beta_i(h^p),
\]
\[
\varphi\beta_i(a) = \sum \varphi(a_h) \varphi(\beta_i h),
\]
\[
ver_i(p\beta_{i-1}(a)) = \sum \varphi(a_h) ver_ip\beta_i(h),
\]
\[
ver_i(\sum \tilde{\omega}_i^k\beta_{i-1}(a)) = \sum \varphi(a_h) ver_i(\sum \tilde{\omega}_i^k\beta_{i-1}(h)).
\]
Hence we need to prove the assertion only for $h$, a conjugacy class of $G/\Gamma^{(e)}$. But this is obvious.   \qed \\

\begin{lemma} For any $u\in \widehat{\lm(\G)_{\ts}}^{\times}$, we have
\[
\frac{\hat{\theta_i}\vp(u)}{\vp\hat{\theta_i}(u)} = \frac{ver_i(\hat{\theta}_{i-1}(u))^p}{\prod_{k=0}^{p-1}ver_i(\tilde{\omega_i^{k}}\hat{\theta}_{i-1}(u))},
\]
for all $ i \geq 1$.
\label{beta2}
\end{lemma}
\noindent{\bf Proof:} $\wh{\lm(\G)_{\ts}}$ may be obtained from $\wh{\lm(\Gi)_{\ts}}$ by adjoining a $p^i$th root of $\g^{p^i}$. $\hat{\theta}_i(u)$ is the norm of $u$ in this extension of rings $\wh{\lm(\G)_{\ts}}/\wh{\lm(\Gi)_{\ts}}$ and so can be given in terms of its ``conjugates". $\vp$ on $\wh{\lm(\Gi)_{\ts}}$ has image in $\wh{\lm(\G^{(i+1)})_{\ts}}$ and coincides with $ver_{i+1}$. Now it can be easily checked that, in fact
\[
\hat{\theta}_i\varphi(u)=ver_i(\hat{\theta}_{i-1}(u))^p, \text{                   and}
\]
\[
\varphi\hat{\theta}_i(u)=\prod_{k=0}^{p-1}ver_i(\tilde{\omega_i^k}\hat{\theta}_{i-1}(u)).
\]
This finishes the proof. \qed \\

\noindent Take any $y \in J(\lm(G))$ (resp. $y \in J(\widehat{\lm(G)_{\ts}})$). Then for any $i \geq 1$:
\begin{align*}
\beta_i(L([1+y]))   &= \beta_i(log([1+y])-\frac{\varphi}{p}log([1+y])) \\
                                &=log(\theta_i([1+y])) - \frac{1}{p}\beta_i(\varphi(log([1+y]))) 
\end{align*}

\begin{align*}
\beta_i(\varphi(log([1+y]))) &= (\varphi(\beta_i)+ver_i(p\beta_{i-1})-ver_i(\sum_{k=0}^{p-1}\tilde{\omega}_i^k(res_{i-1})))(log([1+y]) \\
                                                &= log(\varphi(\theta_i([1+y]))) + log(ver_i(\theta_{i-1}([1+y]))^p)-log(\prod_{k=0}^{p-1}ver_i(\tilde{\omega}_i^k(\theta_{i-1}([1+y]))))  \\
                                                &=log \Big( \frac{\varphi(\theta_i([1+y]))ver_i(\theta_{i-1}([1+y]))^p}{\prod_{k=0}^{p-1}ver_i(\tilde{\omega}_i^k(\theta_{i-1}([1+y])))} \Big).                                                                                                                                              
\end{align*}                                                
Hence
\[
\beta_i(L([1+y])) = \frac{1}{p}log\Big( \Big( \frac{\theta_i([1+y])}{ver_i(\theta_{i-1}([1+y]))}\Big)^p \Big(\frac{\prod_{k=0}^{p-1}ver_i(\tilde{\omega}_i^k(\theta_{i-1}([1+y])))}{\varphi(\theta_i([1+y]))} \Big)\Big).
\]
Similarly, for any $y \in J(\widehat{\lm(G)_{\ts}})$, we have:
\[
\hat{\beta}_i(L([1+y])) = \frac{1}{p}log\Big( \Big( \frac{\hat{\theta}_i([1+y])}{ver_i(\hat{\theta}_{i-1}([1+y]))}\Big)^p \Big(\frac{\prod_{k=0}^{p-1}ver_i(\tilde{\omega}_i^k(\hat{\theta}_{i-1}([1+y])))}{\varphi(\hat{\theta}_i([1+y]))} \Big)\Big).
\] 

\noindent Using lemma \ref{beta2} above, we get:

\begin{proposition} For any $x\in \lm(G)^{\times}$ (resp. $x \in \widehat{\lm(G)_{\ts}}^{\times}$) and any $i\geq 1$ we have:
\begin{equation}
\beta_i(L([x]))= \frac{1}{p}log \Big( \Big(\frac{\theta_i([x])}{ver_i(\theta_{i-1}([x]))}\Big)^p\Big(\frac{ver(\prod_{k=0}^{p-1}\tilde{\omega}_i^k(\theta_{i-1}([x])))}{\varphi(\theta_i([x]))}\Big)\Big)
\label{equation1}
\end{equation}
\begin{equation}
(resp. \ \ \hat{\beta}_i(L([x]))= \frac{1}{p}log \Big( \Big(\frac{\hat{\theta}_i([x])}{ver_i(\hat{\theta}_{i-1}([x]))}\Big)^p\Big(\frac{ver(\prod_{k=0}^{p-1}\tilde{\omega}_i^k(\hat{\theta}_{i-1}([x])))}{\varphi(\hat{\theta}_i([x]))}\Big)\Big))
\label{equation2}
\end{equation}
\end{proposition}

\section{Computation of $K_1$}
In this section we state and prove our main theorems about $K_1$-groups. We continue to use the notations of previous sections.
\subsection{An additive theorem}
We have homomorphisms
\[
\beta=(\beta_i) : T(\lm(G)) \ra \prod_{i=0}^e \lm(G_i^{ab}), \text{       and}
\]
\[
\hat{\beta} : (\hat{\beta}_i) : T(\wh{\lm(G)_{\ts}}) \ra \prod_{i=0}^e\wh{\lm(G_i^{ab})_{\ts}}.
\]
In this section we prove that these homomorphisms are injective and describe their images.

\begin{definition} For $0 \leq j \leq i \leq e$, $H_j\times \Gi$ is a subgroup of $G_j^{ab}$ of finite index. Thus $\lm(G_j^{ab})$ (resp. $\lm(G_j^{ab})_{\ts}$) is a free $\lm(H_j\times \Gi)$-module (resp. $\lm(H_j\times \Gi)_{\ts}$-module) of finite rank. Hence we have the trace map 
\begin{center}
$tr : \lm(G_j^{ab}) \ra \lm(H_j\times \Gi)$\\
$(\text{resp.      } tr: \wh{\lm(G_j^{ab})_{\ts}} \ra \wh{\lm(H_j\times \Gi)_{\ts}}).$
\end{center} 
\end{definition}

\begin{definition} For $0 \leq j \leq i \leq e$ the surjection $G_i^{ab} \ra H_j\times \Gi$ induces the natural map (which is a ring homomorphism)
\begin{center}
$\pi : \lm(G_i^{ab}) \ra \lm(H_j\times \Gi)$\\
$(\text{resp.      } \pi: \wh{\lm(G_i^{ab})_{\ts}} \ra \wh{\lm(H_j\times \Gi)_{\ts}}).$  
\end{center}
\end{definition} 

\begin{definition} $\G$ acts by conjugation on $G_i^{ab}$ and so on $\lm(G_i^{ab})$ and $\lm(G_i^{ab})_{\ts}$, for all $0 \leq i \leq e$. Let $T_i$ (resp. $\hat{T}_i$) be the image of the map on $\lm(G_i^{ab})$ (resp. $\wh{\lm(G_i^{ab})_{\ts}}$) which maps $x$ to $\sum_{k=0}^{p^i-1}\g^kx\g^{-k}$. We note that $T_i$ (resp. $\hat{T}_i$) is an ideal in the ring $\lm(G_i^{ab})^{\G}$ (resp. $\wh{\lm(G_i^{ab})_{\ts}^{\G}}$).
\end{definition}  

\begin{definition} Let $\Psi$ (resp. $\hat{\Psi}$) be the subgroup of $\prod_{i=0}^e\lm(G_i^{ab})$ (resp. $\prod_{i=0}^e\wh{\lm(G_i^{ab})_{\ts}}$) consisting of all tuples $(a_i)$ satisfying the following conditions: \\
A1) We require that $tr(a_j) = \pi(a_i)$ in $\lm(H_j\times \Gi)$ (resp. $\lm(H_j\times \Gi)_{\ts}$) for all $ 0 \leq j \leq i \leq e$. \\
A2) We require that $a_i$ lies in $T_i$ (resp. $\hat{T}_i$) for all $0\leq i \leq e$.
\end{definition}

\begin{lemma} $\beta(T(\lm(G))) \subset \Psi$ and $\hat{\beta}(T(\wh{\lm(G)_{\ts}})) \subset \hat{\Psi}$.
\label{imageofbeta}
\end{lemma}
\noindent{\bf Proof:} Take any $h\g^a$ in $G$, with $h$ in $H$ and $ 0 \leq a \leq p^e-1$. Let $p^l$ be the highest power of $p$ dividing $a$. Then $\beta_i$ and $\hat{\beta}_i$ both map the class of $h\g^a$ to 0 for all $i > l$. For any $ 0 \leq i \leq l$, let $k_i$ be the smallest integer such that $\g^{p^{k_i}}$ fixes the class of $h$ in $H_i$. Then $\beta_i$ and $\hat{\beta}_i$ both map the conjugacy class of $h\g^a$ to the class of 
\[
p^{i-k_i}(\sum_{t=0}^{p^{k_i}-1}(\g^th\g^{-t})\g^a).
\]
A2) follows easily from this. Note that $tr$ sends the class of $h\g^a$ to 0 if $a$ is not divisible by $p^i$ and if $a$ is divisible by $p^i$, then $tr$ sends the class of $h\g^a$ to $p^{i-j}h\g^a$ (in both cases). If $ 0 \leq j \leq i \leq e$, then 
\begin{align*}
\pi(\beta_i([h\g^a])) &= \pi(\sum_{t=0}^{p^i-1}[\g^th\g^{-t}\g^a]) \\
                               &=p^{i-j}(\sum_{t=0}^{p^j-1}[\g^th\g^{-t}\g^a]) \\
                               &=p^{i-j}\beta_j([h\g^a]) \\
                               &=tr(\beta_j([h\g^a])).
\end{align*}
The same holds for $\hat{\beta}_i$ \qed \\

\noindent Next, we define a $\lm(\Ge)$-linear (resp. $\wh{\lm(\Ge)_{\ts}}$-linear) map $\tau: \Psi \ra T(\lm(G))$ (resp. $\hat{\tau} : \hat{\Psi} \ra T(\wh{\lm(G)_{\ts}})$) as follows: Take any $[h\g^a]$ in $G_i^{ab}$. Then $\tau$ (resp. $\hat{\tau}$) maps $[h\g^a]$ to $\frac{1}{p^i}[h\g^a]$ in $\mQ_p\otimes_{\mZ_p}T(\lm(G))$ (resp. in $\mQ_p\otimes_{\mZ_p}T(\wh{\lm(G)_{\ts}})$) if $a\mZ_p = p^i\mZ_p$ and to 0, otherwise (here $a \in \mZ/p^e\mZ$ and when $a$ is divisible by $p^e$, we take $a\mZ_p$ to be $p^e\mZ_p$). Note that $\tau$ and $\hat{\tau}$ are well defined maps (i.e. their images lie in $T(\lm(G))$ and $T(\wh{\lm(G)_{\ts}})$ respectively) by lemma \ref{conjugatesofg} and A2). The following lemma is immediate from the definitions.

\begin{lemma} $\tau \circ \beta$ is identity on $T(\lm(G))$. $\hat{\tau}\circ \hat{\beta}$ is identity on $T(\wh{\lm(G)_{\ts}})$. 
\label{taubeta}
\end{lemma}
\qed 

\begin{lemma} $\tau$ and $\hat{\tau}$ are injective.
\label{tauisoneone}
\end{lemma} 
\noindent{\bf Proof:} Let $(a_i)$ be an element in the kernel of $\tau$. First observe that, when $a_i$ is written (by lemma \ref{crossedproduct}) as a linear combination of elements of $G/\Ge$ with coefficients in $\lm(\Ge)$, the coefficients of $[h\g^a]$ is 0 for all $a$ such that $a\mZ_p = p^i\mZ_p$ (by lemma \ref{conjugatesofg} and definition of $\tau$). Now the trace of $[h\g^a] \in \lm(G_j^{ab})$ in $\lm(H_i\times \G^{(j)})$, with $a\mZ_p = p^i\mZ_p$ is $p^{i-j}[h\g^a]$. So coefficient of $[h\g^a]$ in $a_j$ is 0 by A1). Hence $\tau$ is injective. Similarly, we show that $\hat{\tau}$ is injective.
\qed\\ 
                                
\noindent Putting everything together we get the main theorem of this section.

\begin{theorem} $\beta$ is an isomorphism between $T(\lm(G))$ and $\Psi$. $\hat{\beta}$ is an isomorphism between $T(\wh{\lm(G)_{\ts}})$ and $\hat{\Psi}$. 
\label{additivetheorem}
\end{theorem}
\qed 

\subsection{Description of $K_1$-groups}
\begin{definition} For $0 \leq j \leq i \leq e$, we have the norm map
\begin{center}
$Nr: \lm(G_j^{ab})^{\times} \ra \lm(H_j\times\Gi)^{\times}$ \\
$(\text{resp.    } Nr: \wh{\lm(G_j^{ab})_{\ts}}^{\times} \ra \wh{\lm(H_j\times\Gi)_{\ts}}^{\times}),$
\end{center}
\end{definition} 

\begin{definition} Let $D_i$ (resp. $\hat{D}_i$) be the image of the map on $\lm(G_i^{ab})$ (resp. $\wh{\lm(G_i^{ab})_{\ts}}$) given by $x \mapsto \sum_{k=0}^{p-1}\g^{kp^{i-1}}x\g^{-kp^{i-1}}$.
\end{definition} 

\begin{definition} Let $\Phi$ (resp. $\hat{\Psi}$) be the set of all tuples $(x_i)$ in $\prod_{i=0}^e\lm(G_i^{ab})^{\times}$ (resp. $\prod_{i=0}^e\wh{\lm(G_i^{ab})_{\ts}}^{\times}$) satisfying the following conditions: \\
M1) For $0 \leq j \leq i \leq e$, we require $Nr(x_j)=\pi(x_i)$. \\
M2) $\G$ acts on $H_i$ and hence on $\lm(G_i^{ab})$ and $\wh{\lm(G_i^{ab})_{\ts}}$. We want $x_i$ to be fixed under the action of $\g$.\\
M3) For all $1 \leq i \leq e$, we require
\begin{center}
$x_i \equiv ver_i(x_{i-1}) \ (mod \ D_i+(p))$ \\
$(\text{resp.    } x_i \equiv ver_i(x_{i-1}) \ (mod \ \hat{D}_i+(p))).$
\end{center}
M4) We want 
\[
\Big(\Big(\frac{x_i}{ver_i(x_{i-1})}\Big)^p\Big(\frac{ver_i(\prod_{k=0}^{p-1}\tilde{\omega}_i^k(x_{i-1}))}{\varphi(x_i)}\Big)\Big) \equiv 1 \ (mod \ pT_i)
\]
\[
(\text{resp.    } \Big(\Big(\frac{x_i}{ver_i(x_{i-1})}\Big)^p\Big(\frac{ver_i(\prod_{k=0}^{p-1}\tilde{\omega}_i^k(x_{i-1}))}{\varphi(x_i)}\Big)\Big) \equiv 1 \ (mod \ p\hat{T}_i)).
\]
\label{definitionofphi}
\end{definition}

\noindent It is easily seen that $\Phi$ and $\hat{\Phi}$ are subgroups of $\prod_{i=0}^e\lm(G_i^{ab})^{\times}$ and $\prod_{i=0}^e\wh{\lm(G_i^{ab})_{\ts}}^{\times}$ respectively.

\begin{lemma} We have homomorphisms 
\begin{center}
$L: \Phi \ra \Psi$,     and       $L: \hat{\Phi} \ra \hat{\Psi}$,
\end{center}
given by $(x_i) \mapsto (a_i)$, where
\[
a_i=
\begin{cases}
\frac{1}{p}log\Big(\Big(\frac{x_i}{ver_i(x_{i-1})}\Big)^p\Big(\frac{ver_i(\prod_{k=0}^{p-1}\tilde{\omega}_i^k(x_{i-1}))}{\varphi(x_i)}\Big)\Big) \ \ \ \ \text{if} \  i\geq 1 \\
\frac{1}{p}log\Big(\frac{x_0^p}{\varphi(x_0)}\Big) \hspace{3.7cm} \text{if} \ i=0
\end{cases}
\]
Note that $log$ is defined on the expression in brackets since some power of it is congruent to 1 modulo $p$.
\label{phiintopsi}
\end{lemma} 
\noindent{\bf Proof:} Since $tr \circ log = log \circ Nr$ and $\pi \circ log \cong log \circ \pi$, $(a_i)$ satisfy A1). Since $log$ is an isomorphism between $1+pT_i$ and $pT_i$ (resp between $1+p\hat{T}_i$ and $p\hat{T}_i$), $(a_i)$ satisfy A2).  \qed 

\begin{lemma} $\theta(K_1(\lm(G))) \subset \Phi$ and $\hat{\theta}(K_1(\wh{\lm(G)_{\ts}})) \subset \hat{\Phi}$.
\label{imageoftheta}
\end{lemma} 
\noindent{\bf Proof:} The map $\theta_i$ can be described as follows. $\lm(G)$ is a free $\lm(G_i)$-module with a basis $\{1, \gamma, \gamma^2, \ldots, \gamma^{p^i-1}\}$. Norm of any $x \in \lm(G)^{\times}$ in $K_1(\lm(G_i))$ is the class of the matrix of the $\lm(G_i)$-linear map on $\lm(G)$ given by multiplication by $x$ on the right. Call this matrix $A$. $\theta_i$ just takes $x$ to the determinant of image  of $A$ under $\lm(G_i) \ra \lm(G_i^{ab})$. $\hat{\theta}_i$ has similar description. From this it is easy to see that M1) and M2) are satisfied by the image of $\theta$ and $\hat{\theta}$. M4), i.e. the congruences follow very easily from the relations (1) and (2), and the and definition of $\Psi$ and $\hat{\Psi}$. We now prove M3). We prove it only for $\lm(G)$. The case of $\wh{\lm(G)_{\ts}}$ is similar. For $i \geq 1$, we have following commutative diagram: 

\xymatrix{K_1(\lm(G)) \ar[rd]^{norm} \ar[d]_{norm} & \\
K_1(\lm(G_{i-1})) \ar[d] \ar[r]_{norm} & K_1(\lm(G_i)) \ar[d] \\
\lm(G_{i-1}^{ab})^{\times} & \lm(G_i^{ab})^{\times}}

\noindent It suffices to prove that the square commutes modulo $D_i+(p)$ if we put arrow 
\[
ver_i: \lm(G_{i-1}^{ab})^{\times} \ra \lm(G_i^{ab})^{\times}.
\]
Let $\sigma$ denote the automorphism of $\lm(G_i)$ given by $x \ra \gamma^{p^{i-1}}x\gamma^{-p^{i-1}}$. Let $x = \sum_{k=0}^{p-1} \gamma^{p^{i-1}k}x_k \in \lm(G_{i-1})^{\times}$, with all $x_k \in \lm(G_i)$. Then norm of $[x] \in K_1(\lm(G_{i-1}))$ in $K_1(\lm(G_i))$ is class of following matrix:
\[
\left( \begin{array}{cccc}
x_0 & x_1 & \cdots & x_{p-1} \\
\sigma(x_{p-1}) & \sigma(x_0) & \cdots & \sigma(x_{p-2}) \\
\vdots & & & \vdots \\
\sigma^{p-1}(x_1) & \sigma^{p-1}(x_2) & \cdots & \sigma^{p-1}(x_0)
\end{array} \right)
\] 

\noindent $\theta_i([x])$ is just the determinant of this matrix in $\lm(G_i^{ab})^{\times}$. It can be proved from this that image of $[x]$, say $y$ satisfies 
\[
y \equiv \sum_{k=0}^{p-1}\prod_{i=0}^{p-1}\sigma^i([x_k]) (mod \ D_i +(p)).
\]
And it is easy to see that 
\[
\prod_{i=0}^{p-1}\sigma^i([x_k]) \equiv ver_i([x_k]) (mod D_i + (p)).
\]
This proves the lemma. \qed 

\begin{lemma} There is a short exact sequence
\[
1 \ra G^{ab} \times \mu_{p-1} \ra \Phi \xrightarrow{L} \Psi \ra G^{ab} \ra 1.
\]
Here the first map is $x \mapsto (v_i(x))$, where $v_i$ is the map induced by the transfer homomorphism $G^{ab} \ra G_i^{ab}$, and the last map  takes $a_0$ to an element of $G^{ab}$ given by the exact sequence of lemma \ref{exactsequenceofintegrallog}
\end{lemma}
\noindent{\bf Proof:} Let $(x_i)$ be an element in the kernel of $L$. By the exact sequence of lemma \ref{exactsequenceofintegrallog} we get that $x_i \in \mu_{p-1} \times G_i^{ab}$. It is clear that $(y_i)\in \Phi$, where 
\[
y_i=
\begin{cases}
x_0 \ \ \ \ if \ i=0; \\
v_i(x_0) \ \ if \ i\geq 1.
\end{cases}
\]
here $v_i:\mu_{p-1}\times G^{ab} \ra \mu_{p-1}\times G_i^{ab}$ is the identity map in the first component and the transfer map in the second component. So after replacing $x_i$ by $x_iy_i^{-1}$, we assume that $x_0=1$. Inductively assume that $x_i=1$ for all $0 \leq i\leq k$. Then by M4) we have $L(x_{k+1})=0$. Which means $x_{k+1} \in \mu_{p-1}\times G_{k+1}^{ab}$. Now M3) says that
\[
ver_{k+1}(x_k)=1\equiv x_{k+1} (mod \ D_{k+1}+(p)).
\]
But coefficient of 1 in $D_{k+1}+(p)$ is always a multiple of $p$. Hence $x_i=1$ for every $i$. So an element in the kernel of $L$ is determined by $x_0$. Which proves that exactness at $\Phi$. 

\noindent We now prove exactness at $\Psi$. We have a commutative diagram with exact rows:

\xymatrix{1 \ar[r] & \mu_{p-1} \times G^{ab} \ar@{=}[d] \ar[r] & K_1(\lm(G)) \ar[d]_{\theta_0} \ar[r]^L & T(\lm(G)) \ar[d]^{\beta_0} \ar[r] & G^{ab} \ar@{=}[d] \ar[r] & 1 \\
1 \ar[r] & \mu_{p-1} \times G^{ab} \ar[r] & \lm(G^{ab})^{\times} \ar[r]_L & \lm(G^{ab}) \ar[r]_{\eta} & G^{ab} \ar[r] & 1}

\noindent The map from $\Psi \ra G^{ab}$ is $(a_i) \mapsto \eta(a_0)$. Hence, if some element $(a_i) \in \Psi$, then there exists a $a \in T(\lm(G))$ such that $\beta(a) = (a_i)$. If $(a_i)$ lies in the kernel of $\eta$. Then $a$ must be image of some $x$ under $L$. Then $L(\theta(x)) = (a_i)$. This proves the exactness at $\Psi$.  \qed \\

\begin{theorem} $\theta: K_1(\lm(G)) \ra \Phi$ is an isomorphism.
\label{thetaisisomorphism}
\end{theorem}
\noindent{\bf Proof:} We have a commutative diagram with exact rows:

\xymatrix{1 \ar[r] & G^{ab}\times \mu_{p-1} \ar[r] \ar@{=}[d] & K_1(\lm(G)) \ar[d]_{\theta} \ar[r]^{L} & T(\lm(G)) \ar[d]^{\beta} \ar[r] & G^{ab} \ar[r] \ar@{=}[d] & 1 \\
1 \ar[r] & G^{ab}\times \mu_{p-1} \ar[r] & \Phi \ar[r]_{L} & \Psi \ar[r] & G^{ab} \ar[r] & 1}

\noindent $\beta$ is an isomorphism. A trivial diagram chase shows that $\theta$ is an isomorphism. \qed \\

\subsection{The homomorphism $\theta_{\ts}$ and the group $\Phi_{\ts}$}
Just like $\theta$, we define a homomorphism $\theta_{\ts}=(\theta_{i,\ts}):K_1(\lm(G)_{\ts}) \ra \prod_{i=0}^e\lm(G_i^{ab})_{\ts}^{\times}$. Let $\Phi_{\ts}=\prod_{i=0}^e\lm(G_i^{ab})_{\ts}^{\times}\cap\hat{\Phi}$. Then $\Phi_{\ts}$ can also be described by 4 conditions M1), M2), M3), and M4) as above. Note that they did not involve any use of logarithm. 

\begin{theorem} $\theta_{\ts}(K_1(\lm(G)_{\ts})) \subset \Phi_{\ts}$. and $\Phi_{\ts}\cap \lm(G_i^{ab})^{\times} =\Phi$.
\label{thetat}
\end{theorem}
\noindent{\bf Proof:} Clear. \qed \\

\noindent I am not sure if the following corollary is known.

\begin{corollary} The natural map $K_1(\lm(G))\ra K_1(\lm(G)_{\ts})$ is an injection.
\label{injection}
\end{corollary} 
\noindent{\bf Proof:} We have
\xymatrix{K_1(\lm(G)) \ar[r] \ar[d]^{\cong} & K_1(\lm(G)_T) \ar[d] \\
\Phi \ar@{^{(}->}[r] & \Phi_T}

\noindent Therefore $K_1(\lm(G)) \ra K_1(\lm(G)_T)$ is injective.
\qed  

\subsection{Application} The theorems proven in this section provide us with the algebraic ingredient for proving the Main Conjecture. Define $I$ by 
\[
I=\{ (H\rtimes\Gi, (\g^{p^i}-1)H) : 0 \leq i \leq e\}.
\]
Lemma \ref{representationsofg} shows that C1) is satisfied. Theorem \ref{thetaisisomorphism} gives C2) and theorem \ref{thetat} gives C3).

\section{A Special Case} \label{sectionspecialcase} In this section we prove the noncommutative Main Conjecture in a special case. We continue to use the notation introduced at the beginning of section 2. 

\subsection{Simplification of congruences} 
\begin{lemma} The map $\vp$ from $G_i$ to $G_{i+1}$ induces a homomorphism from $G_i^{ab}$ to $G_{i+1}^{ab}$ if and only if
\[
[h^p] = [h\g^{p^i}h\g^{-p^i} \cdots \g^{(p-1)p^i}h\g^{-(p-1)p^i}]
\]
in $H_{i+i}$, for every $h$ in $H$. In the case when $\vp$ induces a homomorphism from $G_i^{ab}$ to $G_{i+1}^{ab}$ it coincides with the transfer homomorphism $ver_{i+1}$ from $G_i^{ab}$ to $G_{i+1}^{ab}$. 
\label{veri}
\end{lemma} 
\noindent{ \bf Proof:} Take any $g = h\g^a$ in $G_i$. If $p^{i+1}$ divides $a$ then $\vp(g)$ is equal to $\vp(h)\vp(\g^a)$ in $G_{i+1}^{ab}$. So lets assume that $a\mZ_p = p^i\mZ_p$. Then
\[
\vp(g) = (h\g^ah\g^{-a}\cdots \g^{(p-1)a}h\g^{-(p-1)a})\g^{ap}, \text{    and}
\]
\[
\vp(h)\vp(\g^a) = h^p\g^{ap}.
\]
So $\vp$ induces a homomorphism from $G_i^{ab}$ to $G_{i+1}^{ab}$ if and only if $[h^p] = [(h\g^ah\g^{-a}\cdots \g^{(p-1)a}h\g^{-(p-1)a})]$ in $G_{i+1}^{ab}$. But $[(h\g^ah\g^{-a}\cdots \g^{(p-1)a}h\g^{-(p-1)a})]$ is the same as $[(h\g^{p^i}h\g^{-p^i} \cdots \g^{(p-1)p^i}h\g^{-(p-1)p^i})]$ in $G_{i+1}^{ab}$ by lemma \ref{relationing}. 
\[
ver_{i+1}(h\g^a) 
\begin{cases} = (h\g^a)^p \hspace{2cm} \ \   if  \ \  a\mZ_p = p^i \mZ_p \\
                      =(h\g^{p^i}h\g^{-p^i} \cdots \g^{(p-1)p^i}h\g^{-(p-1)p^i})\g^{ap} \ \ \ if \ \  a\mZ_p \neq p^i\mZ_p.
\end{cases}
\]
Therefore, if $\vp$ induces a homomorphism from $G_i^{ab}$ to $G_{i+1}^{ab}$, then it must be the same as $ver_{i+1}$.  \qed 

\begin{definition} We call $G$ of special type if $\vp$ induces a homomorphism from $G_i^{ab}$ to $G_{i+1}^{ab}$ for every integer $i \geq 0$. 
\label{specialtype}
\end{definition}

\begin{lemma} Assume that $G$ is of special type, then
\begin{center}
$ \beta_i(\vp(x)) = pver_i(\beta_{i-1}(x)) $  for all $x \in T(\lm(G))$, and \\
$ \hat{\beta}_i(\vp(x)) = pver_i(\hat{\beta}_{i-1}(x)) $  for all $x \in T(\wh{\lm(G)_{\ts}})$.
\end{center}
\label{beta3}
\end{lemma}
\noindent{\bf Proof:} As in lemma \ref{beta1} we need to check this only when $x$ is a conjugacy class of $G/\Ge$. This is again obvious. \qed \\

\begin{corollary} We have the relations:
\[
\beta_i(L(x)) = log(\frac{\theta_i(x)}{ver_i(\theta_{i-1}(x))}) \ \ for \ all \ x \in K_1(\lm(G)),
\]
\[
\hat{\beta}_i(L(x)) = log(\frac{\hat{\theta}_i(x)}{ver_i(\hat{\theta}_{i-1}(x))}) \ \ for \ all \ x \in K_1(\wh{\lm(G)_{\ts}}).
\]
\label{specialcase}
\end{corollary}
\qed 

\noindent Thus, if $G$ is of special type, we get the following simple description of $\Phi$ and $\hat{\Phi}$. $\Phi$ (resp. $\hat{\Phi}$) consists of all tuples $(x_i)$ in $\prod_{i=0}^e\lm(G_i^{ab})^{\times}$ (resp. $\wh{\lm(G_i^{ab})_{\ts}}^{\times}$) satisfying \\
MS1) $Nr(x_j)=\pi(x_i)$ for all $0 \leq j \leq i \leq e$. \\
MS2) $x_i$ is fixed under the action of $\g$. \\
MS3) $x_i\equiv ver_i(x_{i-1}) \ (mod \ T_i)$ (resp. $(mod \ \hat{T}_i)$). \\
All unexplained notation is as before and we get theorem \ref{thetaisisomorphism} and theorem \ref{thetat}. 

\subsection{Proof of the Main Conjecture in a special case} In this section we prove congruences predicted in the special case considered above. Lets fix the notation. $F$ is a totally real number field of degree $r$ over $\mQ$. $F_{\infty}/F$ is an admissible $p$-adic Lie extension with Galois group $G=H\rtimes \Gamma$, with $H$ a finite abelian $p$-group. We assume that $G$ is of special type. As usual we assume that $F_{\infty}/F$ satisfies hypothesis $\mu =0$. Let $F_0 =F$ and let $F_i$ be the unique extension of $F$ in $F^{cyc}$ of degree $p^n$. Let $K_i$ be the maximal abelian extension of $F_i$ contained in $F_{\infty}$. Hence $K_e = F_{\infty}$ and $Gal(K_i/F_i)= G_i^{ab}$. Let $\Sigma$ be a fixed finite set of primes of $F$ which contains all primes which ramify in $F_{\infty}$ and all primes above $p$. The set of primes in $F_i$ lying above $\Sigma$ will be denoted by $\Sigma_i$. Let $\zeta(K_i/F_i)$ be the $p$-adic abelian zeta function of Deligne and Ribet for the abelian extension $K_i/F_i$. $\g$ acts on $H_i$ and so also on $\\lm(G_i^{ab})_{\ts}$. Let $T_{i,\ts}$ be the image of the map on $\lm(G_i^{ab})_{\ts}$ given by $x \mapsto \sum_{k=0}^{p^i-1}\g^k x \g^{-k}$. It is an ideal in the ring $\lm(G_i^{ab})_{\ts}^{\G}$. We call this ideal the \emph{trace ideal}. 

\begin{theorem} With the notation as above $\zeta(K_i/F_i) \equiv ver(\zeta(K_{i-1}/F_{i-1}))$ modulo the trace ideal. This proves the Main Conjecture for the extension $F_{\infty}/F$. 
\label{congruences}
\end{theorem}

\noindent We will prove this congruence in this section using techniques of Kato \cite{Kato:2006} and of Ritter and Weiss \cite{RitterWeiss:congruences}. Our proof is obtained by a slight extension of techniques used in \cite{RitterWeiss:congruences}. 

\subsubsection{Approximations to $\zeta(K_i/F_i)$} In this section we get a sequence of elements in certain group rings which essentially approximate $\zeta(K_i/F_i)$. These group rings are obtained as follows. Let $N$ be the composition 
\[
Gal(K_0(\mu_p)/F) \twoheadrightarrow Gal(F(\mu_{p^{\infty}}/F)\xrightarrow{\kappa_F} \mZ_p^{\times}.
\]
Then for any positive integer $k$, divisible by $p-1$, $N^{k}$ factors through the group $G^{ab}$. Define $f$ to be the positive integer such that $N^{p-1}(G^{ab})$ is equal to the subgroup $1+p^f\mZ_p$ of $\mZ_p^{\times}$. Then we have an isomorphism
\[
\lm(G_i^{ab}) \xrightarrow{\sim} \ilim{j \geq i}(\mZ_p[G_i^{ab}/\G^{(j)}]/(p^{f+j})).
\]
Since $f+j \ra \infty$, injection is clear. For the surjection, observe that given any $(x_j)_j$ in the projective limit $\varprojlim (\mZ_p[G_i^{ab}/\G^{(j)}]/(p^{f+j}))$, we can canonically construct $\tilde{x}_j$ in $\mZ_p[G_i^{ab}/\G^{(j)}]$ as follows. Let $\bar{x}_t$ be the image of $x_t$ in $\mZ_p[G_i^{ab}/\G^{(j)}]/(p^{f+t})$. Define $\tilde{x}_j$ to be the limit of $\bar{x}_t$, for all $t \geq j$. It is easy to see that $\tilde{x}_j$ from an inverse system to give an element in $\lm(G_i^{ab})$. This element maps to $(x_j)_j$. Hence the above map is also surjective. \\
Let $x$ be a coset of an open subgroup $U$ of $G_i^{ab}$. Set $\delta^{(x)}(g)$ to be 1 if $g \in x$, and to be 0 otherwise. Define the \emph{partial zeta function}, $\zeta_i(\delta^{(x)}, s)$, by 
\[
\zeta_i(\delta^{(x)}, s) = \sum_{\mathfrak{a}}\frac{\delta^{(x)}(g_{\mathfrak{a}})}{Norm(\mathfrak{a})^s}, \hspace{2cm} Re(s) >1,
\]
where the sum is over all ideals $\mathfrak{a}$ of $O_{F_i}$, which are prime to $\Sigma_i$, $g_{\mathfrak{a}}$ is the Artin symbol of the ideal $\mathfrak{a}$ in $G_i^{ab}$, and $Norm$ is the absolute norm of the ideal $\mathfrak{a}$. It is well known from the work of Seigel, that $\zeta_i(\delta^{(x)}, s)$ has analytic continuation to the whole complex plane except for a simple pole at $s=1$, and that $\zeta_i(\delta^{(x)}, 1-k)$ in rational for any even positive integer $k$. If $\epsilon$ is a locally constant function on $G_i^{ab}$ with values in a $\mQ$-vector space $V$, say for some open subgroup $U$ of $G_i^{ab}$, 
\[
\epsilon = \sum_{x \in G_i^{ab}/U} \epsilon(x) \delta^{(x)},
\]
then the $L$-value of $\epsilon$ at $1-k$ can be canonically defined as 
\[
L_i(\epsilon, 1-k) = \sum_x \epsilon(x) \zeta_i(\delta^{(x)}, 1-k) \in V.
\]
If $\epsilon$ is an Artin character of degree 1, then this is of course the value at $1-k$ of the complex $L$-function associated to $\epsilon$ with Euler factors at primes in $\Sigma_i$ removed. If $\epsilon$ is a locally constant $\mQ_p$- valued function on $G_i^{ab}$, then, for any positive integer $k$ divisible by $p-1$, define
\[
\Delta_i(\epsilon, 1-k) = L_i(\epsilon, 1-k) - N(\gamma^{p^i})^kL_i(\epsilon_{(i)}, 1-k),
\]
where, $\epsilon_{(i)}$ is a locally constant function defined by $\epsilon_{(i)}(h)=\epsilon(\gamma^{p^i}h)$, for all $h \in G_i^{ab}$. Recall that $\gamma$ is the fixed topological generator of $\G$. The following proposition is an easy consequences of the congruences of Coates \cite{Coates:1977} which are proven in Deligne-Ribet (\cite{DeligneRibet:1980}, theorems 0.4 and 0.5) (see also Ritter and Weiss \cite{RitterWeiss:congruences}).

\begin{proposition} Under the natural map from $\lm(G_i^{ab})$ to $\mZ_p[G_i^{ab}/\G^{(j)}]/(p^{f+j})$, $(1-\g^{p^i})\zeta(K_i/F_i)$ maps to 
\[
\sum_{x \in G_i^{ab}/\G^{(j)}} \Delta_i(\delta^{(x)}, 1-k)N(x)^{-k}x \ (mod \ p^{f+j}),
\]
for any $j \geq i$ and any positive integer $k$ divisible by $p-1$. In particular, we are claiming that the inverse limit is independent of the choice of $k$.
\label{approximation}
\end{proposition}
\qed 

\subsubsection{A sufficient condition} We reduce the congruence between $p$-adic zeta functions to congruences between the elements of group algebra approximating $(1-\gamma^{p^i})\zeta(K_i/F_i)$. First we prove a small lemma.

\begin{lemma} Let $y$ be a coset of $\Gamma^{(j)}$ in $G_i^{ab}$. Then 
\[
\Delta_i(\delta^{(y)}, 1-k) = \Delta_i(\delta^{(y^{\gamma})}, 1-k).
\]
\end{lemma}
\noindent{\bf Proof:} We just need to prove that $\zeta_i(\delta^{(y)}, 1-k) = \zeta_i(\delta^{(y^{\gamma})}, 1-k)$ because of the following
\[
\Delta_i(\delta^{(y)}, 1-k) = \zeta_i(\delta^{(y)}, 1-k) - N(\gamma^{p^i})^k \zeta_i(\delta^{(y)}_{(i)}, 1-k),
\] 
\[
\Delta_i(\delta^{(y^{\gamma})}, 1-k) = \zeta_i(\delta^{(y^{\gamma})}, 1-k) - N(\gamma^{p^i})^k \zeta_i(\delta^{(y^{\gamma})}_{(i)}, 1-k)
\]
and
\[
\delta^{(y)}_{(i)} = \delta^{(y\gamma^{-p^i})},
\]
\[
\delta^{(y^{\gamma})}_{(i)} = \delta^{((y\gamma^{-p^i})^{\gamma})}.
\]
So we must prove that $\zeta_i(\delta^{(y)}, 1-k) = \zeta_i(\delta^{(y^{\gamma})}, 1-k)$. But this is clearly true as $Norm(\mathfrak{a}) = Norm(\mathfrak{a}^{\gamma})$ and image of $\mathfrak{a}$ (under Artin reciprocity map) lies in $y$ if and only if $\mathfrak{a}^{\gamma}$ lies in $y^{\gamma}$. \qed \\

\noindent Let $V_i$ denote the kernel of $ver_i$, the transfer map from $G_{i-1}^{ab}$ to $G_i^{ab}$. It is clear that $V_i$ is contained in $H_{i-1}$. Let $T_{i,j}$ denote the image of the map on $\mZ_p[G_i^{ab}/\G^{(j)}]/(p^{f+j})$ given by 
\[
x \mapsto \sum_{t=0}^{p^i-1} \g^{p^t}x\g^{-p^t}.
\]
We call this the trace ideal of $\mZ_p[G_i^{ab}/\G^{(j)}]/(p^{f+j})$ (though it is actually an ideal only in the ring $(\mZ_p[G_i^{ab}/\G^{(j)}]/(p^{f+j}))^{\G}$). 

\begin{lemma} A sufficient condition for the theorem \ref{congruences} to hold is that, for any positive integer $k$ divisible by $p-1$, 
\[
\Delta_i(\epsilon, 1-k) \equiv \Delta_{i-1}(\epsilon\circ ver_i, 1-pk) \ (mod \ p^{i-j}),
\]
for all locally constant $\mZ_{(p)}$-valued functions $\epsilon$ on $G_i^{ab}$ fixed by $\g^{p^j}$. Here $\g$ acts on $\epsilon$ as follows: $(\g \cdot \epsilon)(h) = \epsilon(h^{\g})$.
\label{sufficientcondition}
\end{lemma}
\noindent{\bf Proof:} The theorem says that 
\[
\zeta(K_i/F_i) \equiv ver_i(\zeta(K_{i-1}/F_{i-1})) \ (mod \ T_{i, \ts}).
\]
This is easily seen to be equivalent to 
\[
(1-\g^{p^i})\zeta(K-i/F_i) \equiv ver_i((1-\g^{p^{i-1}})\zeta(K_{i-1}/F_{i-1})) \ (mod \ T_i).
\]
Thus we compare the images of $(1-\g^{p^i})\zeta(K_i/F_i)$ and $ver_i((1-\g^{p^{i-1}})\zeta(K_{i-1}/F_{i-1}))$ in $\mZ_p[G_i^{ab}/\G^{(l)}]/(p^{f+l-1})$, for every $l \geq i$. These images are \\
($i$) $\sum_{y \in G_i^{ab}/\G^{(l)}} \Delta_i(\delta^{(y)}, 1-k)N(y)^{-k}y \ (mod \ p^{f+l-1})$, and \\
($ii$) $\sum_{x \in G_{i-1}^{ab}/V_i\times \G^{(l-1)}} \Delta_{i-1}(\delta^{(x)}, 1-pk)N(x)^{-pk}ver_i(x) \ (mod \ p^{f+l-1})$.\\
Here $k$ is any positive integer divisible by $p-1$. We use the independence of $k$ in proposition \ref{approximation}. Assume that $y$ is fixed by $\g^{p^j}$ but not by $\g^{p^{j-1}}$. Then $\delta^{(y)}$ is fixed by $\g^{p^j}$. Using this and the previous lemma, we conclude that the $\g$ orbit of $y$ in ($i$) is
\begin{equation}
\sum_{t=0}^{p^j-1}\Delta_i(\delta^{(y)}, 1-k)N(y)^{-k}\g^ty\g^{-t} \ (mod \ p^{f+l-1}).
\label{sum}
\end{equation}
By hypothesis we have
\[
\Delta_i(\delta^{(y)}, 1-k)\equiv \Delta_{i-1}(\delta^{(y)}\circ ver_i, 1-pk) \ (mod \ p^{i-j}).
\]
Now, if $y$ is not in the image of $ver_i$, then $\delta^{(y)}\circ ver_i=0$ and so $\Delta_{i-1}(\delta^{(y)}\circ ver_i, 1-pk)=0$. Then 
\[
\Delta_{i}(\delta^{(y)}, 1-k) \equiv 0 \ (mod \ p^{i-j}),
\] 
and the above sum (\ref{sum}) lies in the trace ideal of $\mZ_p[G_i^{ab}/\G^{(l)}]/(p^{f+l-1})$. On the other hand, if $y = ver_i(x)$, then $x$ in $G_{i-1}^{ab}/V_i\times\G^{(l-1)}$ is uniquely determined by $y$. This $x$ is fixed by $\g^{p^j}$ and not fixed by $\g^{p^{j-1}}$. Hence the  $\g$ orbit of $x$ in ($ii$) is 
\begin{equation}
\sum_{t=0}^{p^j-1}\Delta_{i-1}(\delta^{(x)}, 1-pk)N(x)^{-pk}ver_i(\g^tx\g^{-t}) \ (mod \ p^{f+l-1}).
\label{sum2}
\end{equation}
Now, $N(y)^{-k} = N(ver_i(x))^{-k} = N(x)^{-pk}$, and $\delta^{(y)}\circ ver_i = \delta^{(x)}$. Hence
\[
\Delta_i(\delta^{(y)}, 1-k) \equiv \Delta_{i-1}(\delta^{(x)}, 1-pk) \ (mod \ p^{i-j}), 
\]
and the difference of sums in (\ref{sum}) and (\ref{sum2}) lies in the trace ideal of $\mZ_p[G_i^{ab}/\G^{(l)}]/(p^{f+l-1})$. The inverse image of these trace ideals over all $l$'s gives the trace ideal $T_i$ of $\lm(G_i^{ab})$ and hence finishes the proof. \qed \\

\subsubsection{Proof of the congruences} We now prove the congruences in lemma \ref{sufficientcondition}. We do this using the technique of Ritter and Weiss which depends on the theory of Deligne and Ribet. We recall some basic notions of this theory. Lets denote the degree of $F$ over $\mQ$ by $r$. We recall the notion of Hilbert modular forms over $F$. Let $\mathfrak{H}_F$ denote the Hilbert upper half plane of $F$. Let $\mathfrak{f}$ be an integral ideal of $F$ with all prime factors in $\Sigma$. We put $GL_2^+(F\otimes \mathbb{R})$ for the group of $2\times 2$ matrices with totally positive determinant. It acts on function $f : \mathfrak{H}_F \ra \mathbb{C}$ by
\[
f|_k\left(\begin{array}{cc} a & b \\ c & d \end{array} \right)(\tau) = \N(ad-bc)^{k/2}\N(c\tau +d)^{-k}f\Big(\frac{a\tau+b}{c\tau+d}\Big), 
\]
where $\N: F\otimes \mathbb{C} \ra \mathbb{C}$ is the norm map. Set 
\[
\G_{00}(\mathfrak{f}) = \{ \left( \begin{array}{cc} a & b \\ c & d \end{array} \right) \in SL_2(F) : a,d \in 1+\mathfrak{f}, b \in \mathfrak{D}^{-1}, c \in \mathfrak{f}\mathfrak{D} \},
\]
where $\mathfrak{D}$ is the different of $F$. A \emph{Hilbert modular form} $f$ of weight $k$ on $\G_{00}(\mathfrak{f})$ is a homomorphic function $f : \mathfrak{H}_F \ra \mathbb{C}$ (which we also assume to be homomorphic at $\infty$ if $F =\mQ$) satisfying
\[
f|_kM=f \hspace{2cm} \text{for all} \ M\in \G_{00}(\mathfrak{f}).
\]
The space of all Hilbert modular form of weight $k$ on $\G_{00}(\mathfrak{f})$ is denoted by $M_k(\G_{00}(\mathfrak{f}), \mathbb{C})$. Any $f$ in $M_k(\G_{00}(\mathfrak{f}), \mathbb{C})$ has the standard $q$ expansion 
\[
c(0) + \sum_{\mu} c(\mu) q_F^{\mu},
\]
where $\mu$ runs through all totally positive elements in $O_F$, and $q_F^{\mu} = e^{2\pi i tr_{F/\mQ}(\mu\tau)}$. The inclusion of $F$ into $F_i$ induces maps $\mathfrak{H}_F \xrightarrow{*} \mathfrak{H}_{F_i}$ and $SL_2(F\otimes\mathbb{R}) \xrightarrow{*} SL_2(F_i\otimes\mathbb{R})$.  For a homomorphic function $f : \mathfrak{H}_{F_i} \ra \mathbb{C}$, we define $R_if$ to be the homomorphic function $R_if : \mathfrak{H}_F \ra \mathbb{C}$ such that $R_if(\tau) = f(\tau^{*})$. Then 
\[
(R_if)|_{p^ik}M = R_i(f|_kM^*),
\]
for any $M  \in SL_2(F\otimes\mathbb{R})$. It is well known that if  $f$ lies in the space $M_k(\G_{00}(\mathfrak{f}O_{F_i}), \mathbb{C})$, for an integral ideal $\mathfrak{f}$ of $F$, then $R_if$ lies in the space $M_{p^ik}(\G_{00}(\mathfrak{f}), \mathbb{C})$, and if the standard $q$-expansion of $f$ is
\[
c(0) + \sum_{\nu \in O_{F_i}} c(\nu)q_{F_i}^{\nu},
\]
then the standard $q$-expansion of $R_if$ is 
\[
c(0) +\sum_{\mu \in O_F}c_*(\mu)q_F^{\mu},
\]
with $c_*(\mu) = \sum_{\nu : tr_{F_i/F}(\nu)=\mu}c(\nu)$. Here $\nu$ and $\mu$ are always totally positive. \\

\noindent Let $\mathbb{A}_{F}$ denote the ring of finite adeles of $F$. Then by strong approximation 
\[
SL_2(\mathbb{A}_{F}) = \wh{\G_{00}(\mathfrak{f})}\cdot SL_2(F).
\]
Any $M$ in $SL_2(\mathbb{A}_{F})$ can be written as $M_1M_2$ with $M_1 \in \wh{\G_{00}(\mathfrak{f})}$ and $M_2 \in SL_2(F)$. We define $f|_kM$ to be $f|_kM_2$.  Any $\alpha$ in $\mathbb{A}_{F}^{\times}$ determines a \emph{cusp}. We let 
\[
f|_{\alpha} = f|_k\left( \begin{array}{cc} \alpha & 0 \\ 0 & \alpha^{-1} \end{array} \right).
\]
The $q$-expansion of $f$ at the cusp determined by $\alpha$ is defined to be the standard $q$-expansion of $f|_{\alpha}$. We write it as
\[
c(0, \alpha) + \sum_{\mu} c(\mu, \alpha) q_F^{\mu},
\]
where the sum is restricted to totally positive elements of $F$ which lie in the square of the ideal ``generated" by $\alpha$. Now let $f$ be any element of $M_k(\G_{00}(\mathfrak{f}O_{F_i}), \mathbb{C})$. Then the constant term of $q$-expansion of $R_if$ at the cusp determined by $\alpha$ in $\mathbb{A}_{F}^{\times}$ is equal to the constant term of $q$-expansion of $f$ at the cusp determined by $\alpha^*$. We need the following lemma proven, for example, in Ritter and Weiss (\cite{RitterWeiss:congruences}, lemma 6)

\begin{lemma} Let $\beta \in O_F$ be a totally positive element. Assume that $\mathfrak{f}\subset \beta O_F$. Then there is a Hecke operator $U_{\beta}$ on $M_k(\G_{00}(\mathfrak{f}), \mathbb{C})$ so that if $f \in M_k(\G_{00}(\mathfrak{f}), \mathbb{C})$ has standard $q$-expansion 
\[
c(0) + \sum_{\mu \in O_F, \mu \gg 0}c(\mu)q_F^{\mu}, 
\]
then $f|_kU_{\beta}$ has the standard $q$-expansion 
\[
c(0) + \sum_{\mu \in O_F, \mu \gg 0}c(\mu\beta)q_F^{\mu}. 
\]
\label{heckeoperator}
\end{lemma}
\qed
  
\noindent The following proposition, which attaches Eisenstein series to a locally constant $\mathbb{C}$-values function $\epsilon$ on $G^{ab}$, is proven in Deligne and Ribet \cite{DeligneRibet:1980} (see section 6, or proposition 8 in Ritter-Weiss \cite{RitterWeiss:congruences}). 

\begin{proposition} Let $\epsilon$ be a locally constant $\mathbb{C}$-valued function on $G^{ab}$. Then for every positive integer $k$ divisible by $p-1$ \\
($i$) There is an integral ideal $\mathfrak{f}$ of $F$ with all its prime factors in $\Sigma$, and a Hilbert modular form $G_{k,\epsilon}$ in $M_k(\G_{00}(\mathfrak{f}), \mathbb{C})$ with standard $q$-expansion
\[
2^{-r}L_i(\epsilon, 1-k) + \sum_{\mu}\Big(\sum_{\mathfrak{a}}\epsilon(\mathfrak{a})N(\mathfrak{a})^{k-1}\Big)q_F^{\mu},
\]
where the first sum ranges over all totally positive $\mu$ in $O_F$, and the second sum ranges over all integral ideals $\mathfrak{a}$ of $F$ containing $\mu$ and prime to $\Sigma$. $\epsilon(\mathfrak{a})$ is defined to be $\epsilon(g_{\mathfrak{a}})$, with $g_{\mathfrak{a}}$ begin the Artin symbol of $\mathfrak{a}$ (see Deligen-Ribet \cite{DeligneRibet:1980}, 2..22). \\
($ii$) The $q$-expansion of $G_{k,\epsilon}$ at the cusp determined any $\alpha$ in $\mathbb{A}_{F}^{\times}$ has constant term
\[
N((\alpha))^k2^{-r}L_i(\epsilon_{g}, 1-k),
\]
where $(\alpha)$ is the ideal of $F$ generated by $\alpha$, $g$ is the image of $\alpha$ under the Artin symbol map, and $\epsilon_g$ is the locally constant function given by 
\[
\epsilon_g(h)=\epsilon(gh) \hspace{2cm} \text{for all} \ h \in G^{ab}.
\]
\end{proposition}
\qed 

\noindent After the above discussion, we are now ready to prove the following key proposition which gives the required congruences.

\begin{proposition} Let $\epsilon$ be a locally constant $\mZ_{(p)}$-valued function on $G_i^{ab}$. Then there exists an integral ideal $\mathfrak{f}$ of $F$ such that all its prime factors are in $\Sigma$ and $\mathfrak{f} \subset p^iO_{F}$, so that for any positive integer $k$, divisible by $p-1$,
\[
E:= R_iG_{k,\epsilon}|_{p^ik}U_{p^i} - R_{i-1}G_{pk, \epsilon\circ ver_i}|_{p^ik}U_{p^{i-1}}, 
\]
belongs to $M_{p^ik}(\G_{00}(\mathfrak{f}), \mathbb{C})$. The constant term of the standard $q$-expansion of $E$ is 
\[
2^{-p^ir}L_i(\epsilon, 1-k)-2^{-p^{i-1}r}L_{i-1}(\epsilon\circ ver_i, 1-pk).
\]
If $\g^{p^j}$ fixes $\epsilon$, then all non-constant coefficients of the standard $q$-expansion of $E$ are in $p^{i-j}\mZ_{(p)}$.
\end{proposition}

\noindent{\bf Proof:} Choose a $\mathfrak{f} \subset p^iO_F$ such that all its prime factors are in $\Sigma$ and $G_{pk, \epsilon\circ ver_i} \in M_{pk}(\G_{00}(\mathfrak{f}O_{F_{i-1}}), \mathbb{C})$ and $G_{k,\epsilon} \in M_k(\G_{00}(\mathfrak{f}O_{F_i}), \mathbb{C})$. So that $E$ lies in $M_{p^ik}(\G_{00}(\mathfrak{f}), \mathbb{C})$. The $q$-expansion of $R_iG_{k,\epsilon}|_{p^ik}U_{p^i}$ is 
\[
2^{-p^ir}L_i(\epsilon, 1-k) + \sum_{\mu}\Big(\sum_{(\mb, \eta)}\epsilon(\mb)N(\mb)^{k-1}\big)q_F^{\mu},
\]
where the first sum ranges over all totally positive elements of $O_F$, and the second sum ranges over all $(\mb, \eta)$ such that $\mb$ is an integral ideal of $F_i$, prime to $\Sigma_i$ and $\eta \in \mb$ is a totally positive element such that $tr_{F_i/F}(\eta) = p^i\mu$. The standard $q$-expansion of $R_{i-1}G_{pk, \epsilon\circ ver_i}|_{p^ik}U_{p^{i-1}}$ is 
\[
2^{-p^{i-1}r}L_{i-1}(\epsilon\circ ver_i, 1-pk) + \sum_{\mu}\Big(\sum_{(\ma, \nu)}\epsilon(\ma O_{F_i})N(\ma)^{pk-1}\Big)q_F^{\mu},
\]
where the first sum ranges over all totally positive elements $\mu$ of $O_F$, and the second sum ranges over all $(\ma, \nu)$, with $\ma$ being an integral ideal of $F_{i-1}$ prime to $\Sigma_{i-1}$, and $\nu \in \ma$ is a totally positive element such that $tr_{F_{i-1}/F}(\nu)=p^{i-1}\mu$. We get $\epsilon(\ma O_{F_i})$ in the sum because $\epsilon\circ ver_i(\ma) = \epsilon(\ma O_{F_i})$. So the constant term $E$ is 
\[
2^{-p^ir}L_i(\epsilon, 1-k) - 2^{-p^{i-1}r}L_{i-1}(\epsilon \circ ver_i , 1-pk).
\]
We now prove that the non-constant terms of $E$ lie in $p^{i-j}\mZ_{(p)}$ whenever $\epsilon$ is fixed by $\g^{p^j}$. Let $t$ be the smallest non-negative integer such that $\g^{p^{j+t}}$ fixes $(\mb, \eta)$. if $t=i-j$, then there is no integral ideal $(\ma, \nu)$, with $\ma$ being the integral ideal of $F_{i-1}$ prime to $\Sigma_{i-1}$ and $\nu \in \ma$, such that $\ma O_{F_i} = \mb$ and $\eta = \nu$. Moreover, the $\g^{p^j}$ orbit of $(\mb, \eta)$ is
\[
\sum_{l=0}^{p^{i-j}-1} \epsilon(\mb^{\g^{p^jl}})N(\mb^{\g^{p^jl}})^{k-1} = p^{i-j} \epsilon(\mb)N(\mb)^{k-1}.
\]
On the other hand, if $ t < i-j$, then there exists $(\mc, \delta)$, with $\mc$ being an integral ideal in $F_{j+t}$, prime to $\Sigma_{j+t}$ and $\delta \in \mc$ is totally real such that $\mc O_{F_i} =\mb$ and $\eta = \delta$. We let $\ma = \mc O_{F_{i-1}}$. Then the difference of $\g^{p^j}$ orbit of $(\mb, \delta)$ and of $(\ma, \delta)$ is 
\begin{align*}
\sum_{l=0}^{p^t-1}\Big( \epsilon(\mb^{\g^{p^jl}})N(\mb^{\g^{p^jl}})^{k-1} - \epsilon((\ma O_{F_i})^{\g^{p^j}l})N(\ma^{\g^{p^jl}} )^{pk-1}\Big) 
   &=p^t\Big(\epsilon(\mb)N(\mb)^{k-1} - \epsilon(\mb)N(\ma)^{pk-1}\Big) \\
   &=p^t\epsilon(\mb) \Big( N(\mc)^{p^{i-j-t}(k-1)} - N(\mc)^{p^{i-j-t-1}(pk-1)}\Big) \\
   &\equiv 0 (mod \ p^{i-j}),
\end{align*}
since $N(\mc)^{p^{i-j-t}(k-1)} \equiv N(\mc)^{p^{i-j-t-1}(pk-1)} \ (mod \ p^{i-j-t})$. This completes the proof. \qed \\

\begin{corollary} Let $\epsilon$ be a locally constant $\mZ_{(p)}$-valued function on $G_i^{ab}$ which is fixed by $\g^{p^j}$. Then 
\[
\Delta_i(\epsilon, 1-k) \equiv \Delta_{i-1}(\epsilon \circ ver_i, 1-pk) \ (mod \ p^{i-j}).
\]
\end{corollary}
\noindent{\bf Proof:} Let $E_1$ be the standard $q$-expansion of $E$. Let $\alpha$ be a finite id\`{e}le of $F$ which is mapped to $\g \in G^{ab}$ by the Artin symbol map. Let $E_{\alpha}$ be the $q$-expansion of $E$ at the cusp determined by $\alpha$. Let $E(\alpha) = N(\alpha_p)^{-p^ik}E_{\alpha}$, where $\alpha_p$ is the ``$p^{th}$-component" of $\alpha$. Then by the theory of Deligne and Ribet (\cite{DeligneRibet:1980}, 0.3, section 5), the constant term of $E_1 - E(\alpha)$ lies in $p^{i-j}\mZ_{(p)}$. This constant term is \\

$\Big(2^{-p^ir}L_i(\epsilon, 1-k) - 2^{-p^{i-1}r}L_{i-1}(\epsilon \circ ver_i, 1-pk)\Big)$ \\ 
$-N(\alpha_p)^{-p^ik}N((\alpha))^{-p^ik}\Big(2^{-p^ir}L_i(\epsilon_{(i)}, 1-k) - 2^{-p^{i-1}r}L_{i-1}((\epsilon \circ ver_i)_{(i)}, 1-pk) \Big)$\\ 

$= 2^{-p^ik}\Big(L_i(\epsilon, 1-k)-N(\g^{p^i})^kL_i(\epsilon_{(i)}, 1-k)\Big)$ \\
$-2^{-p^{i-1}k}\Big(L_{i-1}(\epsilon\circ ver_i, 1-pk) - N(\g^{p^{i-1}})^{pk}L_{i-1}((\epsilon\circ ver_i)_{(i)}, 1-pk)\Big)$ \\

$= 2^{-p^ik}\Delta_i(\epsilon, 1-k) - 2^{-p^{i-1}k}\Delta_{i-1}(\epsilon\circ ver_i, 1-pk)$ \\

$\equiv 2^{-p^ik}(\Delta_i(\epsilon, 1-k) - \Delta_{i-1}(\epsilon\circ ver_i, 1-pk))\equiv 0 (mod \ p^{i-j}).$ \\

\noindent This finished the proof of corollary and hence also the proof of theorem \ref{congruences}. \qed\\ 

\noindent Hence we get the following 

\begin{theorem} Let $F_{\infty}/F$ be an admissible $p$-adic Lie extension (see definition \ref{admissible}) satisfying the hypothesis $\mu = 0$. Assume that $H= Gal(F_{\infty}/F^{cyc})$ is a finite abelian $p$-group, so that $G= Gal(F_{\infty}/F) \cong H\rtimes \G$. We assume that $G$ is of special type (see definition \ref{specialtype}). Then the Main Conjecture for the extension $F_{\infty}/F$ is true.
\label{mcinspecialcase}
\end{theorem}

\noindent In section 6 we give some examples of groups of special type. 

\section{Extensions of dimension greater than 1}
In this section we assume that the admissible $p$-adic Lie extension $F_{\infty}/F$ is such that $H = Gal(F_{\infty}/F^{cyc})$ is an abelian compact $p$-adic Lie group which is pro-$p$. So we have an isomorphism $G \cong H\rtimes \G$. If $J$ is any open subgroup of $H$, then it contains a subgroup which is open in $H$ and normal in $G$ (as $J$ is intersection of $H$ with an open subset of $G$ and any open subset of $G$ containing identity contains an open normal subgroup). Hence we obtain an isomorphism 
\[
\lm(G) \xrightarrow{\sim} \ilim{U} \lm(G/U),
\]
where $U$ runs through open subgroups of $H$ which are normal in $G$. In fact we can find a directed system $\{U_i\}_i$ of open subgroups $U_i$ of $H$ which are normal in $G$ such that $\lm(G) \cong \ilim{i}\lm(G/U_i)$. Put $J_i = G/U_i$. Then $J_i$ is a one dimensional compact $p$-adic Lie group. Plainly, $J_i \cong (H/U_i) \rtimes \G$. Let $K_i$ be the Galois extension of $F$, contained in $F_{\infty}$, such that $Gal(K_i/F) \cong J_i$. In other words, $K_i = F_{\infty}^{U_i}$. We now show that validity of the Main Conjecture for each of the one dimensional admissible $p$-adic Lie extensions $K_i/F$ implies the Main Conjecture for $F_{\infty}/F$.

\begin{lemma} The natural map from $K_1(\lm(G))$ to $K_1(\lm(J_i))$ induces an isomorphism
\[
K_1(\lm(G)) \xrightarrow{\sim} \ilim{i}K_1(\lm(J_i)).
\]
\label{inverselimit}
\end{lemma}
\noindent{\bf Proof:} As $\lm(G)$ and $\lm(J_i)$ are local rings, it is well known that $GL_2$ surjects onto $K_1$ and we have the following isomorphisms
\[
GL_2(\lm(G))/E_2(\lm(G)) \xrightarrow{\sim} K_1(\lm(G)), \text{   and}
\]
\[
GL_2(\lm(J_i))/E_2(\lm(J_i)) \xrightarrow{\sim} K_1(\lm(J_i)),
\]
where $E_2$ denotes the subgroup of all $2\times 2$ elementary matrices (see Curtis and Reiner \cite{CurtisReiner:methodsofrepresentationtheory}). Now consider the following commutative diagram 

\xymatrix{ 1 \ar[r] & E_2(\lm(G)) \ar[r] \ar[d] & GL_2(\lm(G)) \ar[r] \ar[d] & K_1(\lm(G)) \ar[r] \ar[d] & 1 \\
1 \ar[r] & \ilim{i}E_2(\lm(J_i)) \ar[r] & \ilim{i} GL_2(\lm(J_i)) \ar[r] & \ilim{i}K_1(\lm(J_i)) \ar[r] & 1}.

\noindent The top row is exact as we have just observed. The inverse system $\{E_2(\lm(J_i))\}_i$ satisfies the Mittag-Leffler condition. In fact, $E_2(\lm(J_l))$ surjects onto $E_2(\lm(J_i))$, for all $ l \geq i$. Hence the bottom row is exact. It is easy to see that the left and the middle arrow are isomorphisms, hence the right is arrow is an isomorphism. \qed \\

\begin{theorem} Assume that the Main Conjecture for the extensions $K_i/F$ is true for all $i$, then the Main Conjecture for the extension $F_{\infty}/F$ is true. 
\label{dimensiongeq1}
\end{theorem}
\noindent{\bf Proof:} Consider the following commutative diagram with exact rows:

\xymatrix{ 0 \ar[r] & K_1(\lm(G)) \ar[r] \ar[d]_{\sim} & K_1(\lm(G)_S) \ar[d] \ar[r]^{\partial} & K_0(\lm(G), \lm(G)_S) \ar[r] \ar[d] & 0 \\
0 \ar[r] & \ilim{i}K_1(\lm(J_i)) \ar[r] & \ilim{i}K_1(\lm(J_i)_S) \ar[r] & \ilim{i}K_0(\lm(J_i), \lm(J_i)_S) &}.

\noindent Here, the exactness at $\varprojlim K_1(\lm(J_i))$ follows from corollary \ref{injection}. The injection of $K_1(\lm(G))$ into $K_1(\lm(G)_S)$ follows from the injection of $\varprojlim K_1(\lm(J_i))$ into $\varprojlim K_1(\lm(J_i)_S)$. It follows from the discussion in the introduction that $[C(K_{i+1}/F)]$ maps to $[C(K_i/F)]$ under the natural map from $K_0(\lm(J_{i+1}),\lm(J_{i+1})_S)$ to $K_0(\lm(J_i),\lm(J_i)_S)$. Hence the uniqueness of the $p$-adic zeta function satisfying the Main Conjecture gives that $\zeta(K_{i+1}/F)$ maps to $\zeta(K_i/F)$ under the natural map from $K_1(\lm(J_{i+1})_S)$ to $K_1(\lm(J_i)_S)$. So
\[
(\zeta(K_i/F))_i \in \ilim{i}K_1(\lm(J_i)_S).
\]
Let $f$ be any element of $K_1(\lm(G)_S)$ such that 
\[
\partial(f) = -[C(F_{\infty}/F)].
\]
Let $(f_i)_i$ be the image of $f$ under the middle vertical arrow. Then $f_i^{-1}\zeta(K_i/F)$ belongs to $K_1(\lm(J_i))$ for each $i$ and 
\[
(f_i^{-1}\zeta(K_i/F))_i \in \ilim{i} K_1(\lm(J_i)).
\]
Let $u$ be the element of $K_1(\lm(G))$ which maps to $(f_i^{-1}\zeta(K_i/F))_i $ under the left vertical arrow. Then it is clear that $uf$ is the $p$-adic zeta function, $\zeta(F_{\infty}/F)$, which we want. \qed \\ 

\noindent Now we let $G$ be the group $H\rtimes\G$, where $H$ is a pro-$p$ compact $p$-adic Lie group and $\G$ is, as usual, isomorphic to the additive group of $p$-adic integers $\mZ_p$. Recall that we have fixed a topological generator $\g$ of $\G$. We extend some definitions of section 2 to this more general setting. Put $G_i$ for the subgroup $H \rtimes \Gi$. 

\begin{definition} We call $G$ a group of special type if the $p$-power map $\vp$ from $G_i$ to $G_{i+1}$ induces a homomorphism from $G_i^{ab}$ to $G_{i+1}^{ab}$, for all $i \geq 0$. 
\end{definition}

\noindent As in section 2, we again let $H_i$ to the the quotient $H/(\g^{p^i}-1)H$. Just as in lemma \ref{veri} we can show that $G$ is of special type if and only if 
\[
[h^p] = [h \g^{p^i}h\g^{-p^i} \g^{2p^i}h\g^{-2p^i} \cdots \g^{(p-1)p^i}h\g^{-(p-1)p^i}],
\]
in $H_{i+1}$, for all $h$ in $H$ and all $i \geq 0$. As above write $G = \ilim{k}J_k$, with $J_k = G/U_k$, with $U_k$ is an open subgroup of $H$, normal in $G$. Then $J_k \cong (H/U_k)\rtimes \G$. Put $J_{k,i} = (H/U_k)\rtimes \Gi$. Hence $J_{k,i}^{ab} = (H/U_k\cdot (\g^{p^i}-1)H) \times \Gi$. If $G$ is of special type, then it is easy to see that each $J_k$ is of special type. Combining theorems \ref{mcinspecialcase} and \ref{dimensiongeq1}, we get

\begin{theorem} Let $F_{\infty}/F$ be an admissible $p$-adic Lie extension satisfying the hypothesis $\mu = 0$ such that $H =Gal(F_{\infty}/F^{cyc})$ is a pro-$p$ compact abelian $p$-adic Lie group, so that $G =Gal(F_{\infty}/F) \cong H\rtimes\G$. We assume that $G$ is of special type. Then the Main Conjecture for the extension $F_{\infty}/F$ is true. 
\label{proofofmc}
\end{theorem}
\qed 

\section{Examples of groups of special type} We now give some examples of groups of special type. By theorem \ref{proofofmc} above we know that the Main Conjecture is true for extensions of these types.  \\
(1) The group $G = \mZ_p\rtimes \G $ is clearly a group of special type. Recall the example from introduction. Let $F = \mQ(\mu_{37})^{+}$, $F_{\infty}$ is the maximal abelian $37$-extension of $F^{cyc}$ unramified outside the unique prime above 37 in $F^{cyc}$. Then $Gal(F^{cyc}/F) \cong \mZ_{37} \rtimes \G$. Hence the Main Conjecture is true for this extension. \\
(2) More generally, if $G =\mZ_p^{r} \rtimes \G$, with $r \geq 1$ and the action of $\g$ on $\mZ_p^r$ is given by a diagonal matrix i.e. some $A \in GL_r(\mZ_p)$ which give action of $\g$ on $\mZ_p^r$ by $\g\cdot (h_1, \ldots, h_r) = A (h_1,\ldots , h_r)^t$ and which is diagonal such that $A$ modulo $p$ is identity matrix. Then $G$ is of special type. \\
(3) Let $G = \mZ_p^2 \rtimes \G$ with $\g$ acting of $\mZ_p^2$ by the matrix $\left(\begin{array}{cc} 1+p & p \\ 0 & 1+p \end{array} \right)$. Then $G$ is of special type. \\
(4) Let $G = \mZ_p^2\rtimes \G$ with $\g$ acting on $\mZ_p^2$ by the matirx $\left(\begin{array}{cc} 1 & 1 \\ 0 & 1 \end{array}\right).$ Indeed, $G$ is isomorphic to the $p$-adic Heisenberg group. This is the example for which Kato \cite{Kato:2006} first proved the Main Conjecture. \\

\noindent It is easy to construct many more examples of groups of special type. However, it is an interesting question to find out exactly which of these groups can occur as Galois group of an adimissible $p$-adic Lie extensions. For instance, assuming Leopoldt's conjecture, the group in (4) cannot occur as a Galois group of an admissible $p$-adic Lie extension. But, its one dimensional quotients may of course occur as Galois groups of admissible $p$-adic Lie extension.

\bibliographystyle{plain}
\bibliography{mybib}
\end{document}